\documentclass{article}
\usepackage{amsmath,amsthm}
\usepackage{verbatim,amsmath,amscd,amssymb,amsthm}
\begin{document}
\font\germ=eufm10
\def\ssl{\hbox{\germ sl}}
\def\slh{\widehat{\ssl_2}}
\def\ge{\hbox{\germ g}}
\def\aaa{@}
\def\aaa{@}
\title{\Large\bf Inductive Construction of Nilpotent Modules of \\
 Quantum Groups at Roots of Unity }
\author{
Yuuki A\textsc{be} 
\thanks
{
e-mail: yu-abe@sophia.ac.jp
}
\\
Department of Mathematics, 
\\
Sophia University
}
\date{}
\maketitle
\begin{abstract}
The purpose of this paper is to prove that 
we can construct all 
finite dimensional irreducible nilpotent modules of type 1 
inductively by using Schnizer homomorphisms 
for quantum algebra at roots of unity of type 
$A_n$, $B_n$, $C_n$, $D_n$ or $G_2$.
\end{abstract}
\maketitle
\renewcommand{\labelenumi}{$($\roman{enumi}$)$}
\renewcommand{\labelenumii}{$(${\rm \alph{enumii}}$)$}
\font\germ=eufm10
\newcommand{\bm}{{\bf m}}
\newcommand{\bU}{{\bf U}}
\newcommand{\cI}{{\mathcal I}}
\newcommand{\cA}{{\mathcal A}}
\newcommand{\cB}{{\mathcal B}}
\newcommand{\cC}{{\mathcal C}}
\newcommand{\cD}{{\mathcal D}}
\newcommand{\cF}{{\mathcal F}}
\newcommand{\cH}{{\mathcal H}}
\newcommand{\cK}{{\mathcal K}}
\newcommand{\cL}{{\mathcal L}}
\newcommand{\cM}{{\mathcal M}}
\newcommand{\cMod}{{\mathcal M}\!{\it od}}
\newcommand{\cN}{{\mathcal N}}
\newcommand{\cO}{{\mathcal O}}
\newcommand{\cS}{{\mathcal S}}
\newcommand{\cV}{{\mathcal V}}
\newcommand{\fra}{\mathfrak a}
\newcommand{\frb}{\mathfrak b}
\newcommand{\frc}{\mathfrak c}
\newcommand{\frd}{\mathfrak d}
\newcommand{\fre}{\mathfrak e}
\newcommand{\frf}{\mathfrak f}
\newcommand{\frg}{\mathfrak g}
\newcommand{\frh}{\mathfrak h}
\newcommand{\fri}{\mathfrak i}
\newcommand{\frj}{\mathfrak j}
\newcommand{\frk}{\mathfrak k}
\newcommand{\frI}{\mathfrak I}
\newcommand{\fm}{\mathfrak m}
\newcommand{\frn}{\mathfrak n}
\newcommand{\frp}{\mathfrak p}
\newcommand{\fq}{\mathfrak q}
\newcommand{\frr}{\mathfrak r}
\newcommand{\frs}{\mathfrak s}
\newcommand{\frt}{\mathfrak t}
\newcommand{\fru}{\mathfrak u}
\newcommand{\frA}{\mathfrak A}
\newcommand{\frB}{\mathfrak B}
\newcommand{\frF}{\mathfrak F}
\newcommand{\frG}{\mathfrak G}
\newcommand{\frH}{\mathfrak H}
\newcommand{\frJ}{\mathfrak J}
\newcommand{\frN}{\mathfrak N}
\newcommand{\frP}{\mathfrak P}
\newcommand{\frT}{\mathfrak T}
\newcommand{\frU}{\mathfrak U}
\newcommand{\frV}{\mathfrak V}
\newcommand{\frX}{\mathfrak X}
\newcommand{\frY}{\mathfrak Y}
\newcommand{\frZ}{\mathfrak Z}
\newcommand{\rA}{\mathrm{A}}
\newcommand{\rC}{\mathrm{C}}
\newcommand{\rd}{\mathrm{d}}
\newcommand{\rB}{\mathrm{B}}
\newcommand{\rD}{\mathrm{D}}
\newcommand{\rE}{\mathrm{E}}
\newcommand{\rH}{\mathrm{H}}
\newcommand{\rK}{\mathrm{K}}
\newcommand{\rL}{\mathrm{L}}
\newcommand{\rM}{\mathrm{M}}
\newcommand{\rN}{\mathrm{N}}
\newcommand{\rR}{\mathrm{R}}
\newcommand{\rT}{\mathrm{T}}
\newcommand{\rZ}{\mathrm{Z}}
\newcommand{\bbA}{\mathbb A}
\newcommand{\bbC}{\mathbb C}
\newcommand{\bbG}{\mathbb G}
\newcommand{\bbF}{\mathbb F}
\newcommand{\bbH}{\mathbb H}
\newcommand{\bbP}{\mathbb P}
\newcommand{\bbN}{\mathbb N}
\newcommand{\bbQ}{\mathbb Q}
\newcommand{\bbR}{\mathbb R}
\newcommand{\bbV}{\mathbb V}
\newcommand{\bbZ}{\mathbb Z}
\newcommand{\adj}{\operatorname{adj}}
\newcommand{\Ad}{\mathrm{Ad}}
\newcommand{\Ann}{\mathrm{Ann}}
\newcommand{\rcris}{\mathrm{cris}}
\newcommand{\ch}{\mathrm{ch}}
\newcommand{\coker}{\mathrm{coker}}
\newcommand{\diag}{\mathrm{diag}}
\newcommand{\Diff}{\mathrm{Diff}}
\newcommand{\Dist}{\mathrm{Dist}}
\newcommand{\rDR}{\mathrm{DR}}
\newcommand{\ev}{\mathrm{ev}}
\newcommand{\Ext}{\mathrm{Ext}}
\newcommand{\cExt}{\mathcal{E}xt}
\newcommand{\fin}{\mathrm{fin}}
\newcommand{\Frac}{\mathrm{Frac}}
\newcommand{\GL}{\mathrm{GL}}
\newcommand{\Hom}{\mathrm{Hom}}
\newcommand{\hd}{\mathrm{hd}}
\newcommand{\rht}{\mathrm{ht}}
\newcommand{\id}{\mathrm{id}}
\newcommand{\im}{\mathrm{im}}
\newcommand{\inc}{\mathrm{inc}}
\newcommand{\ind}{\mathrm{ind}}
\newcommand{\coind}{\mathrm{coind}}
\newcommand{\Lie}{\mathrm{Lie}}
\newcommand{\Max}{\mathrm{Max}}
\newcommand{\mult}{\mathrm{mult}}
\newcommand{\op}{\mathrm{op}}
\newcommand{\ord}{\mathrm{ord}}
\newcommand{\pt}{\mathrm{pt}}
\newcommand{\qt}{\mathrm{qt}}
\newcommand{\rad}{\mathrm{rad}}
\newcommand{\res}{\mathrm{res}}
\newcommand{\rgt}{\mathrm{rgt}}
\newcommand{\rk}{\mathrm{rk}}
\newcommand{\SL}{\mathrm{SL}}
\newcommand{\soc}{\mathrm{soc}}
\newcommand{\Spec}{\mathrm{Spec}}
\newcommand{\St}{\mathrm{St}}
\newcommand{\supp}{\mathrm{supp}}
\newcommand{\Tor}{\mathrm{Tor}}
\newcommand{\Tr}{\mathrm{Tr}}
\newcommand{\wt}{\mathrm{wt}}
\newcommand{\Ab}{\mathbf{Ab}}
\newcommand{\Alg}{\mathbf{Alg}}
\newcommand{\Grp}{\mathbf{Grp}}
\newcommand{\Mod}{\mathbf{Mod}}
\newcommand{\Sch}{\mathbf{Sch}}\newcommand{\bfmod}{{\bf mod}}
\newcommand{\Qc}{\mathbf{Qc}}
\newcommand{\Rng}{\mathbf{Rng}}
\newcommand{\Top}{\mathbf{Top}}
\newcommand{\Var}{\mathbf{Var}}
\newcommand{\gromega}{\langle\omega\rangle}
\newcommand{\lbr}{\begin{bmatrix}}
\newcommand{\rbr}{\end{bmatrix}}
\newcommand{\forb}{\bigcirc\kern-2.8ex \because}
\newcommand{\forbb}{\bigcirc\kern-3.0ex \because}
\newcommand{\forbbb}{\bigcirc\kern-3.1ex \because}
\newcommand{\SpS}{spectral sequence}
\newcommand\C{\mathbb C}
\newcommand\hh{{\hat{H}}}
\newcommand\eh{{\hat{E}}}
\newcommand\F{\mathbb F}
\newcommand\fh{{\hat{F}}}
\newcommand{\End}{\operatorname{End}}
\newcommand{\Stab}{\operatorname{Stab}}
\newcommand{\mo}{\operatorname{mod}}
\newcommand\pf{\noindent {\bf Proof:  }}
\def\al{\alpha}
\def\bf0{\textbf{0}}
\def\c{\textbf{C}}
\def\cq{\textbf{C}(q)}
\def\d{\delta}
\def\deg{\textrm{deg}}
\def\e{\varepsilon}
\def\ep{\epsilon}
\def\g{\frg}
\def\it{\textit}
\def\l{\lambda}
\def\Lnil{L_{\e}^{\rm{nil}}}
\def\no{\nonumber}
\def\ot{\otimes}
\def\rm{\textrm}
\def\ta{\tilde{a}}
\def\tb{\tilde{b}}
\def\tep{\tilde{\epsilon}}
\def\tg{\tilde{g}}
\def\tm{\tilde{m}}
\def\tV{\widetilde{V}}
\def\tv{\tilde{v}}
\def\tx{\tilde{x}}
\def\ty{\tilde{y}}
\def\tz{\tilde{z}}
\def\ua{U_{\varepsilon}(A_n)}
\def\ub{U_{\varepsilon}(B_n)}
\def\uc{U_{\varepsilon}(C_n)}
\def\ud{U_{\varepsilon}(D_n)}
\def\ue{U_{\varepsilon}}
\def\uf{U_{\varepsilon} ^{\rm fin}}
\def\uq{U_q(\mathfrak g)}
\def\ur{U_{\varepsilon} ^{\rm res}}
\def\q{\quad}
\def\qq{\qquad}
\def\vep{\varepsilon}
\def\vp{\varphi}
\def\w{{\mathcal W}}
\def\z{\textbf{Z}}
\section{Introduction}

Let $\uq$ be the quantum algebra of a finite dimensional 
complex simple Lie algebra $\g$ over $\bbC$.
The theory of $\uq$-modules 
is almost same as the one of $\g$ if $q$ is not a root of unity. 
But, if $q$ is a root of unity, 
it is quite different from the one of $\g$. 
For example, there are the following differences.
\begin{itemize}
\item Finite dimensional modules are not always semisimple.
\item Finite dimensional irreducible modules are not necessarily 
highest or lowest weight modules. 
\item Among finite dimensional irreducible modules, 
there exist maximum dimensional modules.
\end{itemize} 
The theory of $\uq$-modules at roots of unity 
is introduced in \cite{DK}. 

Let $\e$ be a primitive $l$-th root of unity.
The completely classification 
of finite dimensional irreducible $\ue(\g)$-modules is not given yet. 
But, the classification of finite dimensional irreducible 
\it{nilpotent} $\ue(\g)$-modules of type $1$ is already given 
by Lusztig in \cite{L1}, \cite{L2} (see \S 3). 
In particular, it is known that 
these modules are classified by highest weights. 

In \cite{N}, Nakashima discover that we can construct these modules 
by using the modules introduced in  \cite{DJMM} 
if $\g$ is type $A_n$. 
Moreover, in \cite{AN}, we discover that we can construct these modules 
by using the \it{Schnizer modules} introduced in \cite{S1} 
if $\g$ is type $B_n$, $C_n$ or $D_n$. 

In this paper, we construct these modules inductively 
in the case that $\g$ is type $A_n$, $B_n$, $C_n$, $D_n$ or $G_2$ 
by using the Schnizer homomorphisms 
introduced in \cite{S2}. 
Then we can construct finite dimensional irreducible 
nilpotent $\ue(\g)$-modules of type $1$ with highest weight 
$(0, \cdots, 0, \l_k, \cdots, \l_n)$ 
as a submodule of a $l^{(N_n-N_{k-1})}$-dimensional $\ue(\g)$-module, 
where $N_n$ is the number of the positive roots of $\g$ 
and $n$ is the rank of $\g$. 
In particular, these results cover the ones of \cite{AN}. 

The organization of this paper is as follows. 
In \S 2, we review the quantum  algebras at roots of unity.
In \S 3, we introduce the nilpotent modules and 
their classification theorem. 
In \S 4, we introduce the Schnizer homomorphisms. 
Finally, in \S5-\S7, we give inductive construction of all 
finite dimensional irreducible nilpotent $\ue(\g)$-modules of type $1$ 
in the case of $\g=A_n$, $B_n$, $C_n$, $D_n$ or $G_2$. 
\section{Quantum algebras at roots of unity}

\setcounter{equation}{0}
\renewcommand{\theequation}{\thesection.\arabic{equation}}

We fix the following notations. 
Let $\g$ be a finite dimensional simple Lie algebra over $\bbC$
 of type $A_n, B_n, C_n, D_n$ or $G_2$.  
We set $I:=\{1,2, \cdots, n\}$. 
Let $ \{\alpha_i\}_{i \in I}$ be the set of simple roots, 
$\Delta$ be the set of roots 
and $\Delta_+$ be the set of positive roots of $\g$.
Let $N$ be the number of positive roots of $\g$, that is, 
$N=\frac{1}{2}n(n+1)$ (resp. $n^2$, $n^2$, $(n-1)n$, $6$) 
if $\g=A_n$ (resp. $B_n$, $C_n$, $D_n$, $G_2$).
We define the root lattice $Q:=\bigoplus_{i \in I} \bbZ \alpha_i$
and the positive loot lattice 
$Q_+:=\bigoplus_{i \in I} \bbZ_+\alpha_i$, 
where $\bbZ_{+}:=\{0, 1, 2, \cdots\}$. 
 Let $(\texttt{a}_{i,j})_{i,j \in I}$ be 
the Cartan matrix associated with $\g$ such that 
\begin{eqnarray*}
  \begin{cases}
   \texttt{a}_{1,2}=-2 & \rm{$\g=B_n$}, \\
   \texttt{a}_{2,1}=-2 & \rm{$\g=C_n$}, \\
   \texttt{a}_{1,2}=0, \, \texttt{a}_{1,3}=\texttt{a}_{2,3}=-1 &
   \rm{$\g=D_n$}, \\
   \texttt{a}_{1,2}=-3 & \rm{$\g=G_2$}. 
  \end{cases}
\end{eqnarray*}
We define $(d_1, \cdots , d_n):=(1, \cdots, 1)$ 
(resp. $(\frac{1}{2}, 1, \cdots, 1)$, 
$(2, 1, \cdots, 1)$, $(1, \cdots, 1)$, $(1, 3)$) 
if $\g=A_n$ (resp. $B_n$, $C_n$, $D_n$, $G_2$). 
We denote the Weyl group of $\g$ by  $\w$ which is generated by the
 simple reflections  $\{s_i\}_{i \in I}$.  

Let $l$ be an odd integer which is greater than $2$. 
We assume that $l$ is not divisible by $3$ if $\g=G_2$. 
Let $\e$ (resp. $\e^{\frac{1}{2}}$) 
be a primitive $l$-th root of unity 
if $\g \neq B_n$ (resp. $\g=B_n$). 
For $r  \in \bbZ$, $m \in \bbN$, $d \in \bbQ$ 
such that $\e^{2d} \neq 1$, we define 
\begin{eqnarray*}
&& [r]_{\e^d}:=\displaystyle \frac{\e^{dr}-\e^{-dr}}{\e^{d}-\e^{-d}}, 
\q [r]:=[r]_{\e}, \\
&& [m]_{\e^{d}}!:=[m]_{\e^{d}}[m-1]_{\e^{d}} \cdots [1]_{\e^{d}}, 
\q  [0]_{\e^{d}}!:=1.
\end{eqnarray*}
\newtheorem{def QA}{Definition}[section]
\begin{def QA}
\label{def QA}
The quantum algebra $\ue(\g)$ is an associative 
$\bbC$-algebra generated by \\
$\{e_{i}, f_{i},t_{i}^{\pm 1} \}_{i \in I}$ 
with the relations    
\begin{eqnarray*}
&& t_{i} t_{i}^{-1} =t_{i}^{-1} t_{i}=1, \q t_{i} t_{j} =
t_{j} t_{i}, \\ 
&& t_{i} e_{j} t_{i}^{-1}=\e_i^{\texttt{a}_{i,j}}e_{j},  
\q  t_{i} f_{j} t_{i}^{-1}=\e_i^{-\texttt{a}_{i,j}}f_{j}, \\ 
&& e_{i} f_{j}- f_{j} e_{i} = \delta_{i,j} \{t_{i}\}_{\e_i}, \\
&& \sum_{k=0}^{1-\texttt{a}_{ij}} (-1)^k e_{i}^{(k)} e_{j} 
e_{i}^{(1-\texttt{a}_{i,j}-k)}=
 \sum_{k=0}^{1-\texttt{a}_{i,j}} (-1)^k f_{i}^{(k)} f_{j} 
f_{i}^{(1-\texttt{a}_{i,j}-k)}=0
\q i \neq j, 
\end{eqnarray*}
where
\begin{eqnarray*}
 e_{i}^{(k)}:= \displaystyle \frac{1}{[k]_{\e_i}!} e_{i}^k, 
\q f_{i}^{(k)} := \displaystyle \frac{1}{[k]_{\e_i}!} f_{i}^k, 
\q \{t_{i}\}_{\e_i} := \displaystyle
\frac{t_i-t_i^{-1}}{\e_i-\e_i^{-1}}, 
\q \e_i:=\e^{d_i}.
\end{eqnarray*}
\end{def QA}
Let $\ue^{+}(\g)$  (resp. $\ue^{-}(\g)$, $\ue^{0}(\g)$) 
be the $\bbC$-subalgebra of $\ue(\g)$ generated by 
$\{e_{i}\}_{i \in I}$ 
(resp. $\{f_{i}\}_{i \in I}$,$\{t_{i}^{\pm 1}\}_{i \in I}$). 
Moreover, we extend the algebra by adding the elements 
$\{t_{i}^{\pm\frac{1}{k}} \, | \, i, k \in I\}$. 

Let $w_0$ be a longest element of $\w$ 
and $w_0= s_{i_1} \cdots s_{i_{N}}$ be a reduced expression of $w_0$. 
We set 
\begin{eqnarray*}
 \beta_1:=\al _{i_1}, \,  \beta_2 :=s_{i_1}(\al_{i_2}), \, 
\cdots , \, \beta_N:=s_{i_1}\cdots s_{i_{N-1}} (\al_{i_{N}}). 
\end{eqnarray*}
Indeed, we have $\Delta_+=\{\beta_1, \cdots , \beta_{N}\}$. 
Then there exist the vectors $\{e_{\beta_i}\}_{i=1}^{N}$, 
$\{f_{\beta_i}\}_{i=1}^{N}$ in $\ue(\g)$ 
which are called ``root vectors'' 
(cf. \cite{DK}, \cite{L2}), 
where $e_{\al_i}=e_{i}$, $f_{\al_i}=f_{i}$ for $i \in I$.
These vectors satisfy the following properties. 
\newtheorem{thm PBW}[def QA]{Proposition}
\begin{thm PBW}[\cite{DK} Proposition 1.7, \cite{L2}]
\label{thm PBW}
\begin{enumerate}
\item
$\{e_{\beta _1}^{m_{1}} \cdots e_{\beta _{N}}^{m_{N}} | 
 m_1, \cdots ,m_{N} \in \bbZ _+ \}$ is a
$\bbC$-basis of $U_{\e}^+ (\g)$.  
\item
$\{f_{\beta_1}^{m_1} \cdots f_{\beta_{N}}^{m_{N}}  | 
 m_1, \cdots , m_{N} \in \bbZ _+ \}$ is a
$\bbC$-basis of $U_{\e}^-(\g)$.  
\item
$\{t_{1}^{m_1} \cdots t_{n}^{m_n}  | 
m_1, \cdots , m_n \in \bbZ\}$
 is a $\bbC$-basis of $U_{\e}^0(\g).$ 
\item
Let $\phi: U_{\e}^-(\g) \otimes U_{\e}^0(\g) \otimes
 U_{\e}^+(\g)\longrightarrow \ue(\g)$  
($u_- \otimes u_0 \otimes u_+   \mapsto u_- u_0  u_+$) 
be the multiplication map. 
Then $\phi$ is an isomorphism of $\bbC$-vector space.
\end{enumerate}
\end{thm PBW}
Let $Z(\ue(\g))$ be the center of $\ue(\g)$.
\newtheorem{pro CE}[def QA]{Proposition}
\begin{pro CE}[\cite{DK} Corollary 3.1]
\label{pro CE}
We have $\{e_{\al}^l , f_{\al}^l , t_{i}^l \, |
\, \al \in \Delta _+,  i \in I\} \subset Z(\ue(\g))$.
\end{pro CE}
Now, for $i \in I$, we set 
\begin{eqnarray}
\deg ( e_{i}):=\al_i, \q \deg ( f_{i}):=-\al_i, 
\q \deg ( t_{i}):=0.
\label{deg}
\end{eqnarray}
Obviously, these are compatible with the relations of $\ue(\g)$.
Therefore, we can regard $\ue(\g)$ as 
$Q$-graded algebra, and we have
\[ \ue(\g)=\bigoplus_{\al \in Q}\ue(\g)_{\al},\q
\ue(\g)_{\al}\ue(\g)_{\al ^{'}} 
\subset \ue(\g)_{(\al + \al ^{'})},\] 
for $\al, \al ^{'} \in Q$,
where $\ue(\g)_{\al}:=\{u \in \ue | \deg (u)= \al \}$. 
\newtheorem{pro deg}[def QA]{Proposition}
\begin{pro deg}[\cite{J} \S8]
\label{pro deg} 
We have $e_{\al} \in U_{\e}^+(\g) \cap \ue(\g)_{\al}$ 
and $f_{\al} \in U_{\e}^-(\g) \cap \ue(\g)_{- \al}$ 
for all $\al \in \Delta_+$. 
\end{pro deg}
\section{Nilpotent modules}
\setcounter{equation}{0}
\renewcommand{\theequation}{\thesection.\arabic{equation}}
\newtheorem{def NM}{Definition}[section]
\begin{def NM}
\label{def NM}
Let $L$ be a $\ue(\g)$-module. 
If $e_{\al}^l=f_{\al}^l=0$ on $L$ for all $\al \in \Delta_{+}$, 
then we call $L$ ``nilpotent module''. 
In particular, if $t_{i}^l=1$ on $L$ for all $i \in I$, 
then we call $L$ ``nilpotent module of type 1''. 
\end{def NM}
\newtheorem{rem type}[def NM]{Remark} 
\begin{rem type}
\label{rem type}
Nilpotent $\ue(\g)$-modules of type 1 
are same as $\uf(\g)$-modules of type 1, 
where $\uf(\g)$ is the finite dimensional quantum algebra  
introduced in \cite{L1}, \cite{L2} (see \cite{AN}). 

In general, finite dimensional irreducible $\uf(\g)$-modules
are divided into $2^n$ types according to 
\{$\sigma : Q \longrightarrow \{\pm 1\}$; homomorphism of group \}.
Without a loss of generality, we may assume that  
finite dimensional irreducible $\uf(\g)$-modules are of type $1$.
\end{rem type}
\newtheorem{def HWM}[def NM]{Definition}
\begin{def HWM}
\label{def HWM}
Let $L$ be a $\ue(\g)$-module.
\begin{enumerate}
\item
We set 
$P(L):=\{v \in L \, | \, e_{i}v=0 \, \, 
\textrm{for all}\, \, i \in I\}$ 
and  call the vectors in $P(L)$ ``primitive vector''. 
\item 
Let $\l=(\l_i)_{i \in I} \in \bbC^n$. 
We assume that $L$ is generated by a nonzero vector $v_0 \in P(L)$ 
such that $t_{i}v_0=\e_i^{\l_i}v_0$ for all $i \in I$. 
Then we call $L$ ``highest weight module with highest weight $\l$'' 
and $v_0$ ``highest weight vector''. 
\end{enumerate}
\end{def HWM}
Now, we introduce the classification theorem of 
finite dimensional
irreducible nilpotent $\ue(\g)$-modules of type 1. 
We set $\bbZ_l:=\{\l \in \bbZ \, | \, 0 \leq \l \leq l-1\}$.
\newtheorem{thm CT}[def NM]{Theorem}
\begin{thm CT}[\cite{L1}, \cite{L2}]
\label{thm CT}
For any $\l \in  \bbZ_l^n$, 
there exists a unique (up to isomorphic) 
finite dimensional irreducible nilpotent $\ue(\g)$-module $\Lnil (\l)$
of type 1 with highest weight $\l$.   
Conversely, if $L$ is a finite dimensional
irreducible nilpotent $\ue(\g)$-module of type 1, 
then there exists a $\l \in  \bbZ_l^n$ 
such that $L$ is isomorphic to $\Lnil (\l)$.
\end{thm CT} 
By the similar manner to the proof of 
Theorem 5.5(ii) in \cite{N} or Theorem 4.10 in \cite{AN}, 
we obtain the following proposition. 
\newtheorem{pro AN}[def NM]{Proposition}
\begin{pro AN}
\label{pro AN}
For $\l \in \bbZ_l^n$, 
let $L$ be a nilpotent highest weight $\ue(\g)$-module of type 1
with highest weight $\l$.
We assume $\textrm{dim}(P(L))=1$. 
Then $L$ is irreducible $\ue(\g)$-module. 
In particular, 
$L$ is isomorphic to $\Lnil (\l)$ as $\ue(\g)$-module.
\end{pro AN}
\section{Schnizer homomorphisms}

\setcounter{equation}{0}
\renewcommand{\theequation}{\thesection.\arabic{equation}}

In the rest of the paper, 
we denote $\g$ by $\g_n$ if the rank of $\g$ is $n$
and $e_i, f_i, t_i$ in $\ue(\g_n)$ by $e_{i,n}, f_{i,n}, t_{i,n}$.

We fix the following notations. 
Let $V_n$ be a $l^n$-dimensional $\bbC$-vector spase 
and $\{v_n(m_n) \, | \, m_n=(m_{1,n}, \cdots,  m_{n,n}) \in \bbZ_l^n\}$ 
be a basis of $V_n$, 
where $\bbZ_l:=\{m \in \bbZ \, | \, 0 \leq m \leq l-1\}$.  
We set $v_n(m_n+lm_n^{'}):=v_n(m_n)$ 
for $m_n, m_n^{'} \in \bbZ_l^n$. 
For $ i \in I$, we set
\begin{eqnarray}
\ep_{i,n}:=(\delta_{i,1}, \delta_{i,2}, \cdots , \delta_{i,n}) \in 
\bbZ_l^n, 
\end{eqnarray}
where $\delta_{i,j}$ is the Kronecker's delta. 
For $i \in I$, $a_{i,n} \in \bbC^{\times}$, $b_{i,n} \in \bbC$, 
we define linear maps 
$x_{i,n}, z_{i,n} \in \textrm{End} (V_n)$ by 
\begin{eqnarray} 
x_{i,n}v_n(m_n):=a_{i,n}v_n(m_n-\ep_{i,n}), 
\q z_{i,n}v_n(m_n):=\e^{m_{i,n}+b_{i,n}}v_n(m_n) 
\q (m_n \in \bbZ_l^n). 
\label{def xz}
\end{eqnarray}
For any $z \in \textrm{End}(V_n)$ such that $z^{-1} \in
\textrm{End}(V_n)$ 
and $d \in \bbQ$ such that $\e^{2d} \neq 1$, we set 
\begin{eqnarray}
 \{z\}_{\e^d}:=\frac{z-z^{-1}}{\e^d-\e^{-d}}.
\label{sig zed}
\end{eqnarray}
Then we have 
\begin{eqnarray}
\{z_{i,n}\}_{\e^d}v_n(m_n)=[d^{-1}(m_{i,n}+b_{i,n})]_{\e^d}v_n(m_n). 
\label{def z}
\end{eqnarray}
For any $\bbC$-vector space $V$, we regard 
$\textrm{End}(V) \otimes \ue(\g_n)$ as $\bbC$-algebra by 
\begin{eqnarray*}
(x \otimes u)(x^{'} \otimes u^{'}):=(xx^{'}) \otimes (uu^{'})
\q (x,x^{'} \in \textrm{End}(V), u,u^{'} \in \ue(\g_n)).
\end{eqnarray*}
\newtheorem{thm SH}{Theorem}[section]
\begin{thm SH}[\cite{S2} Theorem 3.2, 4.10]
\label{thm SH}
(a) Let $\l_n \in \bbC$, 
$a_n=(a_{i,n})_{i=1}^n \in (\bbC^{\times})^n$, 
and $b_n=(b_{i,n})_{i=1}^n \in \bbC^n$. 
Then we obtain a $\bbC$-algebra homomorphism 
$\rho_{n}^A:=\rho_{n}^A(a_n,b_n,\l_n) : \ua
\longrightarrow \rm{End}(V_n) \otimes U_{\e}(A_{n-1})$ 
such that
\begin{eqnarray}
&& \rho_{n}^A(e_{i,n})=\{z_{i-1,n}z_{i,n}^{-1}\}x_{i,n}+
x_{i-1,n}^{-1}x_{i,n}e_{i-1,n-1},
\label{thm SHA1} \\
&& \rho_{n}^A(t_{i,n})=z_{i-1,n}z_{i,n}^{-2}z_{i+1,n}t_{i,n-1},
\label{thm SHA2} \\
&& \rho_{n}^A(f_{i,n})=\{z_{i,n}z_{i+1,n}^{-1}t_{i,n-1}^{-1}\}
x_{i,n}^{-1}+f_{i,n-1},\label{thm SHA3}
\end{eqnarray} 
where 
\begin{eqnarray}
 t_{n,n-1}:=\e^{-\l_n}\prod_{i=1}^{n-1}
t_{i,n-1}^{-\frac{i}{n}}. 
\label{thm SHA4}
\end{eqnarray}

(b) Let $\l_n \in \bbC$, 
$a_n=(a_{i,n})_{i=1}^n  \in (\bbC^{\times})^n$, 
$\ta_{n-1}=(\ta_{i,n-1})_{i=1}^{n-1}  \in (\bbC^{\times})^{n-1}$,
$b_n=(b_{i,n})_{i=1}^n  \in \bbC^n$, 
and $\tb_{n-1}=(\tb_{i,n-1})_{i=1}^{n-1}  \in \bbC^{n-1}$. 
Then we obtain a $\bbC$-algebra homomorphism 
$\rho_{n}^B:=\rho_n^B(a_n,\ta_{n-1},b_{n},\tb_{n-1},\l_n) :
 \ub \longrightarrow \rm{End}(V_n) \otimes  \rm{End}(\tV_{n-1}) 
\otimes  U_{\e}(B_{n-1})$ such that
\begin{eqnarray}
  \rho_{n}^B(e_{i,n})
&=&\{z_{i+1,n}z_{i,n}^{-1}\}x_{i,n}+
\{\tz_{i-1,n-1}\tz_{i,n-1}^{-1}\}x_{i+1,n}^{-1}x_{i,n}\tx_{i,n-1}
\no\\
&& +x_{i+1,n}^{-1}x_{i,n}\tx_{i,n-1}\tx_{i-1,n-1}^{-1}e_{i,n-1} 
\q (2 \leq i \leq n), 
\nonumber \\
 \rho_{n}^B(e_{1,n})
&=&\{z_{2,n}z_{1,n}^{-\frac{1}{2}}\}_{\e_1}x_{1,n}
+\{\tz_{1,n-1}^{-1}z_{1,n}^{\frac{1}{2}}\}_{\e_1}
x_{2,n}^{-1}x_{1,n}\tx_{1,n-1}+x_{2,n}^{-1}\tx_{1,n-1}e_{1,n-1}, 
\label{thm SHB1} \\
\rho_{n}^B(t_{i,n})
&=&z_{i+1,n}z_{i,n}^{-2}z_{i-1,n}
\tz_{i-2,n-1}\tz_{i-1,n-1}^{-2}\tz_{i,n-1}
t_{i,n-1} \q (2 \leq i \leq n), \no \\
 \rho_{n}^B(t_{1,n})
&=&z_{2,n}z_{1,n}^{-1}\tz_{1,n-1}t_{1,n-1}, 
\label{thm SHB2} \\
 \rho_{n}^B(f_{i,n})
&=&\{z_{i,n}z_{i-1,n}^{-1}\tz_{i-2,n-1}^{-1}
\tz_{i-1,n-1}^{2}\tz_{i,n-1}^{-1}
t_{i,n-1}^{-1}\}x_{i,n}^{-1}
\nonumber \\
&& +\{\tz_{i-1,n-1}\tz_{i,n-1}^{-1}t_{i,n-1}^{-1}\}\tx_{i-1,n-1}^{-1}
+f_{i,n-1} \q (2 \leq i \leq n),
\nonumber \\
 \rho_{n}^B(f_{1,n})
&=&\{z_{1,n}^{\frac{1}{2}}\tz_{1,n-1}^{-1}
t_{1,n-1}^{-1}\}_{\e_1}x_{1,n}^{-1}+f_{1,n-1},
\label{thm SHB3} 
\end{eqnarray} 
where 
\begin{eqnarray}
 t_{n,n-1}:=\e^{-\l_n}\prod_{i=1}^{n-1}t_{i,n-1}^{-1}. 
\label{thm SHB4}
\end{eqnarray}

(c) Let $\l_n \in \bbC$, 
$a_n=(a_{i,n})_{i=1}^n  \in (\bbC^{\times})^n$, 
$\ta_{n-1}=(\ta_{i,n-1})_{i=1}^{n-1}  \in (\bbC^{\times})^{n-1}$,
$b_n=(b_{i,n})_{i=1}^n  \in \bbC^n$, 
and $\tb_{n-1}=(\tb_{i,n-1})_{i=1}^{n-1}  \in \bbC^{n-1}$. 
Then we obtain a $\bbC$-algebra homomorphism 
$\rho_{n}^C:=\rho_n^C(a_n,\ta_{n-1},b_{n},\tb_{n-1},\l_n) :
 \uc \longrightarrow \rm{End}(V_n) \otimes  \rm{End}(\tV_{n-1}) 
\otimes  U_{\e}(C_{n-1})$ such that 
$\rho_{n}^C(e_{i,n})$, $\rho_{n}^C(t_{i,n})$, $\rho_{n}^C(f_{i,n})$ 
as in (\ref{thm SHB1}), (\ref{thm SHB2}), (\ref{thm SHB3}) 
if $3 \leq i \leq n$, and
\begin{eqnarray}
 \rho_{n}^C(e_{2,n})
&=&\{z_{3,n}z_{2,n}^{-1}\}x_{2,n}
+\{\tz_{1,n-1}\tz_{2,n-1}^{-1}\}x_{3,n}^{-1}x_{2,n}\tx_{2,n-1}
 +x_{3,n}^{-1}x_{2,n}\tx_{2,n-1}\tx_{1,n-1}^{-1}e_{2,n-1}, \no \\
 \rho_{n}^C(e_{1,n})
&=&\{z_{1,n}^2\tz_{1,n-1}^{-2}\}_{\e_1}x_{2,n}^{-2}
x_{1,n}\tx_{1,n-1}^2 
+\{z_{2,n}\tz_{1,n-1}^{-1}\}x_{2,n}^{-1}x_{1,n}\tx_{1,n-1}\no \\
&&+\{z_{2,n}^2z_{1,n}^{-2}\}_{ \e_1}x_{1,n}
+x_{2,n}^{-2}\tx_{1,n-1}^2e_{1,n-1}, 
\label{thm SHC2} \\
 \rho_{n}^C(t_{2,n})
&=&z_{3,n}z_{2,n}^{-2}z_{1,n}^2
\tz_{1,n-1}^{-2}\tz_{2,n-1}t_{i,n-1},
\q  \rho_{n}^C(t_{1,n})
=z_{2,n}^2z_{1,n}^{-4}\tz_{1,n-1}^2t_{1,n-1}, 
\label{thm SHC3} \\
 \rho_{n}^C(f_{2,n})
&=&\{z_{2,n}z_{1,n}^{-2}\tz_{1,n-1}^{2}\tz_{2,n-1}^{-1}
t_{2,n-1}^{-1}\}x_{2,n}^{-1}
 +\{\tz_{1,n-1}\tz_{2,n-1}^{-1}t_{2,n-1}^{-1}\}\tx_{1,n-1}^{-1}
+f_{2,n-1},\no \\
 \rho_{n}^C(f_{1,n})
&=&\{z_{1,n}^{2}\tz_{1,n-1}^{-2}t_{1,n-1}^{-1}\}_{\e_1}
x_{1,n}^{-1}+f_{1,n-1},
\label{thm SHC4}
\end{eqnarray}
\begin{eqnarray}
t_{n,n-1}:=\e^{-\l_n}t_{1,n-1}^{-\frac{1}{2}}
\prod_{i=2}^{n-1}t_{i,n-1}^{-1} \q (n \geq 2), 
\q t_{1,0}:=\e^{-\l_1}. 
\label{thm SHC1} 
\end{eqnarray} 

(d) Let $\l_n \in \bbC$, 
$a_n=(a_{i,n})_{i=1}^n  \in (\bbC^{\times})^n$, 
$\ta_{n-2}=(\ta_{i,n-2})_{i=1}^{n-2}  \in (\bbC^{\times})^{n-2}$,
$b_n=(b_{i,n})_{i=1}^n  \in \bbC^n$, 
and $\tb_{n-2}=(\tb_{i,n-2})_{i=1}^{n-2}  \in \bbC^{n-2}$. 
Then we obtain a $\bbC$-algebra homomorphism 
$\rho_{n}^D:=\rho_n^D(a_n,\ta_{n-2},b_{n},\tb_{n-2},\l_n) :
 \ud \longrightarrow \rm{End}(V_n) \otimes  \rm{End}(\tV_{n-2}) 
\otimes  U_{\e}(D_{n-1})$ such that 
$\rho_{n}^D(e_{i,n})$,
$\rho_{n}^D(t_{i,n})$, $\rho_{n}^D(f_{i,n})$ as in 
(\ref{thm SHB1}), (\ref{thm SHB2}), (\ref{thm SHB3}) 
if $4 \leq i \leq n$ 
(replace $\tx_{i,n-1}$ to $\tx_{i-1,n-2}$ and 
$\tz_{i,n-1}$ to $\tz_{i-1,n-2}$). 
Moreover,  
\begin{eqnarray}
  \rho_{n}^D(e_{3,n})
&=&\{z_{4,n}z_{3,n}^{-1}\}x_{3,n}
+\{\tz_{1,n-2}\tz_{2,n-2}^{-1}\}x_{4,n}^{-1}x_{3,n}\tx_{2,n-2}
 +x_{4,n}^{-1}x_{3,n}\tx_{2,n-2}\tx_{1,n-2}^{-1}e_{3,n-1}, \no \\
 \rho_{n}^D(e_{2,n})
&=&\{z_{3,n}z_{2,n}^{-1}\}x_{2,n} 
+\{z_{1,n}\tz_{1,n-2}^{-1}\}\tx_{1,n-2}x_{3,n}^{-1}x_{2,n}
 +x_{3,n}^{-1}x_{2,n}x_{1,n}^{-1}\tx_{1,n-2}e_{2,n-1}, \no \\
  \rho_{n}^D(e_{1,n})
&=&\{z_{3,n}z_{1,n}^{-1}\}x_{1,n} 
+\{z_{2,n}\tz_{1,n-2}^{-1}\}\tx_{1,n-2}x_{3,n}^{-1}x_{1,n}
+x_{3,n}^{-1}x_{1,n}x_{2,n}^{-1}\tx_{1,n-2}e_{2,n-1}, \no \\
\label{thm SHD2} \\
\rho_{n}^D(t_{3,n})
&=&z_{4,n}z_{3,n}^{-2}z_{2,n}z_{1,n}
\tz_{1,n-2}^{-2}\tz_{2,n-2}t_{i,n-1}, \no \\
 \rho_{n}^D(t_{2,n})
&=&z_{3,n}z_{2,n}^{-2}\tz_{1,n-2}t_{2,n-1},  
\q \rho_{n}^D(t_{1,n})
=z_{3,n}z_{1,n}^{-2}\tz_{1,n-2}t_{1,n-1},
\label{thm SHD3} \\
 \rho_{n}^D(f_{3,n})
&=&\{z_{3,n}z_{2,n}^{-1}z_{1,n}^{-1}
\tz_{1,n-2}^{2}\tz_{2,n-2}^{-1}t_{3,n-1}^{-1}\}x_{3,n}^{-1}
 +\{\tz_{1,n-2}\tz_{2,n-2}^{-1}t_{3,n-1}^{-1}\}\tx_{1,n-2}^{-1}
+f_{3,n-1}, \no \\
 \rho_{n}^D(f_{2,n})
&=&\{z_{2,n}\tz_{1,n-2}^{-1}t_{2,n-1}^{-1}\}x_{2,n}^{-1}+f_{2,n-1},  \no \\
  \rho_{n}^D(f_{1,n})
&=&\{z_{1,n}\tz_{1,n-2}^{-1}t_{1,n-1}^{-1}\}x_{1,n}^{-1}+f_{1,n-1},
\label{thm SHD4}
\end{eqnarray}
\begin{eqnarray}
t_{n,n-1}:=\e^{-\l_n}t_{1,n-1}^{-\frac{1}{2}}
t_{2,n-1}^{-\frac{1}{2}}\prod_{i=3}^{n-1}t_{i,n-1}^{-1}, 
\q (n \geq 3), 
\q t_{2,1}:=\e^{-\l_2}t_{1,1}^{-\frac{1}{2}}, 
\q t_{1,0}:=\e^{-\l_1}. 
\label{thm SHD1} 
\end{eqnarray}

(e) Let $\l_2 \in \bbC$, 
$a_2=(a_{i,2})_{i=1}^5  \in (\bbC^{\times})^5$, 
and $b_2=(b_{i,2})_{i=1}^5  \in \bbC^5$. 
Then we obtain a $\bbC$-algebra homomorphism 
$\rho^G:=\rho^G(a_2, b_{2}, \l_2) :
 \ue(G_2) \longrightarrow \rm{End}(V_5) \otimes  U_{\e}(A_1)$ such that 
\begin{eqnarray}
\rho^G(e_{1,2})&=&
\{z_{3,2}^3z_{4,2}^{-2}\}x_{1,2}^{-1}x_{2,2}^{-1}x_{3,2}x_{4,2}^2 
+\{z_{4,2}z_{5,2}^{-3}\}x_{1,2}^{-1}x_{2,2}^{-1}x_{4,2}^2x_{5,2}
+\{z_{2,2}^2z_{3,2}^{-3}\}x_{1,2}^{-1}x_{2,2}x_{3,2}  \no \\
&&+[2]\{z_{2,2}z_{4,2}^{-1}\}x_{1,2}^{-1}x_{3,2}x_{4,2} 
+\{z_{1,2}^3z_{2,2}^{-1}\}x_{2,2}
+x_{1,2}^{-1}x_{2,2}^{-1}x_{4,2}x_{5,2}e_{1,1}, \no \\  
\rho^G(e_{2,2})&=&
\{z_{1,2}^{-3}\}_{\e^3}x_{1,2}, 
\label{thm SHG2} \\
\rho^G(t_{1,2})&=&
z_{1,2}^3z_{3,2}^3z_{5,2}^3z_{2,2}^{-2}z_{4,2}^{-2}t_{1,1}, 
\q \rho^G(t_{2,2})=
z_{1,2}^{-6}z_{3,2}^{-6}z_{5,2}^{-6}z_{2,2}^3z_{4,2}^3\e^{\l_2}
t_{1,1}^{-\frac{3}{2}}, 
\label{thm SHG3} \\
\rho^G(f_{1,2})&=&
\{z_{2,2}z_{3,2}^{-3}z_{5,2}^{-3}z_{4,2}^2t_{1,1}^{-1}\}x_{2,2}^{-1}
+\{z_{4,2}z_{5,2}^{-3}t_{1,1}^{-1}\}x_{4,2}^{-1}+f_{1,1}, \no \\
\rho^G(f_{2,2})&=&
\{z_{1,2}^3z_{3,2}^6z_{5,2}^6z_{2,2}^{-3}z_{4,2}^{-3}\e^{-\l_2}
t_{1,1}^{\frac{3}{2}}\}_{\e^3}x_{1,2}^{-1}
+\{z_{3,2}^3z_{5,2}^6z_{4,2}^{-3}\e^{-\l_2}
t_{1,1}^{\frac{3}{2}}\}_{\e^3}x_{3,2}^{-1} \no \\
&&+\{z_{5,2}^3\e^{-\l_2}t_{1,1}^{\frac{3}{2}}\}_{\e^3}x_{5,2}^{-1}. 
\label{thm SHG4}
\end{eqnarray}
Here, $x_{i,j}^{\pm 1}:=0$, $z_{i,j}^{\pm 1}:=1$ 
if the index $(i,j)$ is out of range, 
$e_{0,n}:=f_{n,n-1}:=0$, 
$\ue(\g_0):=\bbC$, $V_0:=\bbC$,
and $\tV_{j}$, $\tx_{i,j}$, $\tz_{i,j}$ are a copy of 
$V_j, x_{i,j}, z_{i,j}$. 
\end{thm SH} 
By using these homomorphisms, 
we obtain $l^{N}$-dimensional $\ue(\g_n)$-modules 
having $\rm{dim}\g_n$-parameters. 
We call those modules the Schnizer modules. 
\newtheorem{rem SH}[thm SH]{Remark}
\begin{rem SH}
\label{rem SH}
The actions of $e_{i,n},t_{i,n},f_{i,n}$ in \cite{S2}
 are slightly
different from the one of Theorem \ref{thm SH}. Because we use a
$\ue(\g_n)$-automorphism $\omega$ such that $(\omega(e_{i,n}),
\omega(t_{i,n}),\omega(f_{i,n}))=(f_{i,n},t_{i,n}^{-1},e_{i,n})$.
\end{rem SH} 
Now, we introduce the following fact to use later. 
If $\g_n=A_n$ $(n \geq 2)$, by (\ref{thm SHA2}), (\ref{thm SHA4}), 
we obtain
\begin{eqnarray}
 \rho_{n-1}^A(t_{n,n-1})
&=&\e^{-\l_n}
\prod_{i=1}^{n-1}\rho_{n-1}^A(t_{i,n-1}^{-\frac{i}{n}}) 
=\e^{-\l_n}\prod_{i=1}^{n-1}
(z_{i-1,n-1}z_{i,n-1}^{-2}z_{i+1,n-1}t_{i,n-2})^{-\frac{i}{n}}  
\nonumber \\
&=&\e^{-\l_n}z_{n-1,n-1}\prod_{i=1}^{n-1}t_{i,n-2}^{-\frac{i}{n}} 
 =\e^{-\l_n}z_{n-1,n-1} t_{n-1,n-2}^{-\frac{n-1}{n}} 
\prod_{i=1}^{n-2}t_{i,n-2}^{-\frac{i}{n}} 
\nonumber \\
 &=&\e^{-\l_n}z_{n-1,n-1} 
(\e^{-\l_{n-1}}\prod_{i=1}^{n-2}t_{i,n-2}^{-\frac{i}{n-1}})^{-\frac{n-1}{n}}  
\prod_{i=1}^{n-2}t_{i,n-2}^{-\frac{i}{n}} 
\nonumber \\
&=&\e^{-\l_n+\frac{n-1}{n}\l_{n-1}}z_{n-1,n-1}. 
\label{3.3.1}
\end{eqnarray}
Similarly, by (\ref{thm SHB2}), (\ref{thm SHB4}), 
(\ref{thm SHC3}), (\ref{thm SHC1}), (\ref{thm SHD3}), (\ref{thm SHD1}), 
we obtain 
\begin{eqnarray}
 \rho_{n-1}^B(t_{n,n-1})
&=&\e^{-\l_n+\l_{n-1}}z_{n-1,n-1}\tz_{n-2,n-2} \q (n \geq 2),
\label{fact tB}\\
 \rho_{n-1}^C(t_{n,n-1})
&=&\e^{-\l_n+\l_{n-1}}z_{n-1,n-1}\tz_{n-2,n-2} \q (n \geq 3), \no \\
\q  \rho_{1}^C(t_{2,1})
&=&\e^{-\l_2+\frac{1}{2}\l_{1}}z_{1,1}^2, 
\label{fact tC}\\
 \rho_{n-1}^D(t_{n,n-1})
&=&\e^{-\l_n+\l_{n-1}}z_{n-1,n-1}\tz_{n-3,n-3} \q (n \geq 4),\no \\ 
 \rho_{2}^D(t_{3,2})
&=&\e^{-\l_3+\frac{1}{2}\l_{2}+\frac{1}{4}}
z_{1,1}^{\frac{1}{2}}z_{1,2}\tz_{2,2},  
\q  \rho_{1}^D(t_{2,1})
=\e^{-\l_2+\frac{1}{2}\l_{1}}z_{1,1},
\label{fact tD}
\end{eqnarray}
$\rho_{n-1}^{\g}(t_{1,0})=\e^{-\l_1}$ for $\g=A, B, C$ or $D$. 
\section{Construction of $\Lnil (\l)$ (Type G-case)}

\setcounter{equation}{0}
\renewcommand{\theequation}{\thesection.\arabic{equation}}

In this section, 
we construct all finite dimensional irreducible nilpotent 
$\ue(G_2)$-modules of type 1 
by using the Schnizer-homomorphisms in Theorem \ref{thm SH}(e). 

We set  
\begin{eqnarray}
&&a_{1,1}^{(0)}:= a_{i,2}^{(0)}:=1 \q (1 \leq i \leq 5),\no \\
&&b_{1,1}^{(0)}:=b_{1,2}^{(0)}:=1, \q b_{2,2}^{(0)}:=4, 
\q b_{4,2}^{(0)}:=5, \q b_{3,2}^{(0)}:=3, \q b_{5,2}^{(0)}:=2, \no \\
&& a_2^{(0)}:=(a_{i,2}^{(0)})_{i=1}^5, 
\q b_2^{(0)}:=(b_{i,2}^{(0)})_{i=1}^5.
\label{sign Gab}
\end{eqnarray} 
For $\l \in \bbC$, we set 
\begin{eqnarray}
&&\rho_1^{A}(\l):=\rho_1^A(a_{1,1}^{(0)}, b_{1,1}^{(0)}, \l) 
: U_{\e}(A_1) \longrightarrow \rm{End}(\bbC), \no \\
&&\rho^{G}(\l):=\rho^G(a_2^{(0)}, b_2^{(0)}, \l) 
: U_{\e}(G_2) \longrightarrow \rm{End}(V_5) \otimes U_{\e}(A_1), 
\label{def rhoGll}
\end{eqnarray}
(see  Theorem \ref{thm SH}(a), (e)). 
For $\l_1, \l_2 \in \bbC$, 
we define 
\begin{eqnarray}
\phi_{1,2}:=\phi_{1,2}(\l_1, \l_2):= 
  \rho_{1}^{A_1}(\nu_1^{(\l_1, \l_2)}) \circ  \rho^{G}(\nu_2^{(\l_1,\l_2)}):
 \ue(G_2) \longrightarrow  \textrm{End} (V_5 \otimes V_1), 
\label{sig Gphikn}
\end{eqnarray}
where 
\begin{eqnarray*}
\nu_{1}^{(\l_1, \l_2)}:=\l_1+2, 
\q \nu_2^{(\l_1, \l_2)}:=\frac{3}{2}\l_1+3\l_2+9.
\end{eqnarray*}
We denote the $\ue(G_2)$-modules associated with 
$(\phi_{1,2}(\l_1,\l_2), V_5 \otimes V_1)$ by $V_{1,2}(\l_1,\l_2)$. 
For $m_{1,1} \in \bbZ_l$, 
$m_5=(m_{i,5})_{i=1}^5 \in \bbZ_l^5$, 
we set
\begin{eqnarray*}
&& v_{1,2}(m_5, m_{1,1}):=v_5(m_5) \otimes v_1(m_{1,1}), 
\q v_{1,2}^{\bf0}:=v_{1,2}(0, \cdots , 0)
\in V_{1,2}(\l_1,\l_2).
\end{eqnarray*}
We define  
$y_{i,2}$, $y_{1,1} \in \textrm{End}(V_{1,2}(\l_1, \l_2))$ 
$(1 \leq i \leq 5)$: for $v=v_{1,2}(m_5, m_{1,1})$, 
\begin{eqnarray} 
y_{1,2}v
&:=&[m_{1,2}+2m_{3,2}+2m_{5,2}-m_{2,2}-m_{4,2}-m_{1,1}-\l_2]_{\e_2}
v_{1,2}(m_5+\ep_{1,2},m_{1,1}), \no \\
y_{2,2}v
&:=&[m_{2,2}-3m_{3,2}-3m_{5,2}+2m_{4,2}+2m_{1,1}-\l_1]
v_{1,2}(m_5+\ep_{2,2},m_{1,1}), \no \\
y_{3,2}v
&:=&[m_{3,2}+2m_{5,2}-m_{4,2}-m_{1,1}-\l_2]_{\e_2}
v_{1,2}(m_5+\ep_{3,2},m_{1,1}), \no \\
y_{4,2}v
&:=&[m_{4,2}-3m_{5,2}+2m_{1,1}-\l_1]
v_{1,2}(m_5+\ep_{4,2},m_{1,1}), \no \\
y_{5,2}v
&:=&[m_{5,2}-m_{1,1}-\l_2]_{\e_2}
v_{1,2}(m_5+\ep_{5,2},m_{1,1}), \no \\
y_{1,1}v
&:=&[m_{1,1}-\l_1]v_{1,2}(m_5,m_{1,1}+\ep_{1,1}). 
\label{def SHGy}
\end{eqnarray}
Then, by Theorem \ref{thm SH}(a), (e), 
(\ref{def xz}), (\ref{def z}), (\ref{sign Gab}),  
we have 
\begin{eqnarray}
e_{1,2}v
&=&[3m_{3,2}-2m_{4,2}]
v_{1,2}(m_{1,2}+1, m_{2,2}+1, m_{3,2}-1, m_{4,2}-2, m_{5,2}, m_{1,1}) \no \\
&&+[m_{4,2}-3m_{5,2}]
v_{1,2}(m_{1,2}+1, m_{2,2}+1, m_{3,2}, m_{4,2}-2, m_{5,2}-1, m_{1,1})\no \\
&&+[2m_{2,2}-3m_{3,2}]
v_{1,2}(m_{1,2}+1, m_{2,2}-1, m_{3,2}-1, m_{4,2}, m_{5,2}, m_{1,1}) \no \\
&&+[2][m_{2,2}-m_{4,2}]
v_{1,2}(m_{1,2}+1, m_{2,2}, m_{3,2}-1, m_{4,2}-1, m_{5,2}, m_{1,1}) \no \\
&&+[3m_{1,2}-m_{2,2}]v_{1,2}
(m_{1,2}, m_{2,2}-1, m_{3,2}, m_{4,2}, m_{5,2}, m_{1,1}) \no \\
&&+[-m_{1,1}]v_{1,2}
(m_{1,2}+1, m_{2,2}+1, m_{3,2}, m_{4,2}-1, m_{5,2}-1, m_{1,1}-1),
\label{fact SHGe1} \\
e_{2,2}v
&=&[-m_{1,2}]_{\e_2}v_{1,2}
(m_{1,2}-1, m_{2,2}, m_{3,2}, m_{4,2}, m_{5,2}, m_{1,1}), 
\label{fact SHGe2} \\
t_{1,2}v
&=&\e^{3m_{1,2}+3m_{3,2}+3m_{5,2}-2m_{2,2}-2m_{4,2}-2m_{1,1}+\l_1}
v_{1,2}(m_5,m_{1,1}), 
\label{fact SHGt1} \\
t_{2,2}v
&=&\e_2^{-2m_{1,2}-2m_{3,2}-2m_{5,2}+m_{2,2}+m_{4,2}+m_{1,1}+\l_2}
v_{1,2}(m_5,m_{1,1}), 
\label{fact SHGt2} \\
f_{1,2}v
&=&(y_{2,2}+y_{4,2}+y_{1,1})v_{1,2}(m_5, m_{1,1}),
\label{fact SHGf1} \\
f_{2,2}v
&=&(y_{1,2}+y_{3,2}+y_{5,2})v_{1,2}(m_5, m_{1,1}).
\label{fact SHGf2}
\end{eqnarray} 
Let $P(V_{1,2}(\l_1, \l_2))$ as in Definition \ref{def HWM} (i).
\newtheorem{sec5}{Proposition}[section]
\begin{sec5}
\label{sec5}
For all $\l_1$, $\l_2 \in \bbC$, 
we obtain $P(V_{1,2}(\l_1,\l_2))=\bbC v_{1,2}^{\bf0}$.
\end{sec5}

{\sl Proof.} 
Since the actions of $e_{1,2}$, $e_{2,2}$ on $V_{1,2}(\l_1, \l_2)$ 
do not depend on $\l_1$, $\l_2$, 
we simply denote $V_{1,2}(\l_1,\l_2)$ by $V_{1,2}$. 
By (\ref{fact SHGe1}), (\ref{fact SHGe2}), obviously, 
$\bbC v_{1,2}^{\bf0} \subset P(V_{1,2})$. 
So we shall prove $P(V_{1,2}) \subset \bbC v_{1,2}^{\bf0}$. 
Let  
\begin{eqnarray*}
 v=\sum_{m_5 \in \bbZ_l^5, m_{1,1} \in \bbZ_l}
c(m_5, m_{1,1})v_{1,2}(m_5, m_{1,1})\in V_{1,2},
\end{eqnarray*}
where 
$c(m_5, m_{1,1}) \in \bbC$, 
and we assume that $e_{1,2}v=e_{2,2}v=0$.  
By (\ref{fact SHGe2}), we get
\begin{eqnarray*}
0=e_{2,2}v
=\sum_{m_5 \in \bbZ_l^5, m_{1,1} \in \bbZ_l}
c(m_5, m_{1,1})[-m_{1,2}]v_{1,2}
(m_{1,2}-1, m_{2,2}, m_{3,2}, m_{4,2}, m_{5,2}, m_{1,1}).
\end{eqnarray*}
Hence, we obtain $c(m_5, m_{1,1})=0$ if $m_{1,2} \neq 0$. 
So, by (\ref{fact SHGe1}), we have
\begin{eqnarray}
0=e_{1,2}v
&=&\sum_{m_{2,2}, m_{3,2}, m_{4,2}, m_{5,2}, m_{1,1} \in \bbZ_l}
c(0, m_{2,2}, m_{3,2}, m_{4,2},m_{5,2}, m_{1,1}) \no \\
&&\{[3m_{3,2}-2m_{4,2}]
v_{1,2}(1, m_{2,2}+1, m_{3,2}-1, m_{4,2}-2, m_{5,2}, m_{1,1}) \no \\
&&+[m_{4,2}-3m_{5,2}]
v_{1,2}(1, m_{2,2}+1, m_{3,2}, m_{4,2}-2, m_{5,2}-1, m_{1,1}) \no \\
&&+[2m_{2,2}-3m_{3,2}]
v_{1,2}(1, m_{2,2}-1, m_{3,2}-1, m_{4,2}, m_{5,2}, m_{1,1}) \no \\
&&+[2][m_{2,2}-m_{4,2}]
v_{1,2}(1, m_{2,2}, m_{3,2}-1, m_{4,2}-1, m_{5,2}, m_{1,1}) \no \\
&&+[-m_{2,2}]
v_{1,2}(0, m_{2,2}-1, m_{3,2}, m_{4,2}, m_{5,2}, m_{1,1}) \no \\
&&+[-m_{1,1}]
v_{1,2}(1, m_{2,2}+1, m_{3,2}, m_{4,2}-1, m_{5,2}-1, m_{1,1}-1)\}. \no 
\end{eqnarray}
Since the $(1,2)$-component of 
$(0, m_{2,2}-1, m_{3,2}, m_{4,2}, m_{5,2}, m_{1,1})$ is $0$ 
and the one of other vectors is $1$, 
by the linearly independence, 
$c(m_5, m_{1,1})=0$ if $m_{2,2} \neq 0$. 
Therefore we obtain 
\begin{eqnarray}
0=e_{1,2}v
&=&\sum_{m_{3,2}, m_{4,2}, m_{5,2}, m_{1,1} \in \bbZ_l}
c(0, 0, m_{3,2}, m_{4,2},m_{5,2}, m_{1,1}) \no \\
&&\{[3m_{3,2}-2m_{4,2}]
v_{1,2}(1, 1, m_{3,2}-1, m_{4,2}-2, m_{5,2}, m_{1,1}) \no \\
&&+[m_{4,2}-3m_{5,2}]
v_{1,2}(1, 1, m_{3,2}, m_{4,2}-2, m_{5,2}-1, m_{1,1}) \no \\
&&+[-3m_{3,2}]
v_{1,2}(1, -1, m_{3,2}-1, m_{4,2}, m_{5,2}, m_{1,1}) \no \\
&&+[2][-m_{4,2}]
v_{1,2}(1, 0, m_{3,2}-1, m_{4,2}-1, m_{5,2}, m_{1,1}) \no \\
&&+[-m_{1,1}]
v_{1,2}(1, 1, m_{3,2}, m_{4,2}-1, m_{5,2}-1, m_{1,1}-1)\}. \no 
\end{eqnarray}
Since the $(2,2)$-component of 
$(1, -1, m_{3,2}-1, m_{4,2}, m_{5,2}, m_{1,1})$ 
(resp. $(1, 0, m_{3,2}-1, m_{4,2}-1, m_{5,2}, m_{1,1})$) 
is $-1$ (resp. $0$) 
and the one of other vectors is $1$, 
we get $c(m_5, m_{1,1})=0$ if $m_{3,2} \neq 0$ or $m_{4,2}=0$. 
Hence we have 
\begin{eqnarray}
0=e_{1,2}v
&=&\sum_{m_{5,2}, m_{1,1} \in \bbZ_l} 
c(0, 0, 0, 0, m_{5,2}, m_{1,1}) 
\{[-3m_{5,2}]v_{1,2}(1, 1, 0, -2, m_{5,2}-1, m_{1,1}) \no \\
&&\qq \qq\qq 
+[-m_{1,1}]v_{1,2}(1, 1, 0, -1, m_{5,2}-1, m_{1,1}-1)\}. \no 
\end{eqnarray}
Since the $(4,2)$-component of 
$(1, 1, 0, -2, m_{5,2}-1, m_{1,1})$ is $-2$ 
and the one of 
$(1, 1, 0, -1$, $ m_{5,2}-1, m_{1,1}-1)$ is $-1$,
we obtain
$c(m_5, m_{1,1})=0$ if $m_{5,2} \neq 0$ or $m_{1,1}=0$. 
It amount to
$v=c(0, \cdots , 0)v_{1,2}^{\bf0} \in \bbC v_{1,2}^{\bf0}$.
\qed \\

Let $y_{i,2}$ $(1 \leq i \leq 5)$, $y_{1,1}$ 
as in (\ref{def SHGy}) 
and we set
$Y_{1,2}:=\{y_{i,2}, y_{1,1} \,| \, 1 \leq i \leq 5 \}$. 
Let $p_{\bf0}: V_{1,2}(\l)\longrightarrow \bbC v_{1,2}^{\bf0}$ 
be the projection.  \\
\newtheorem{lem Gsyaei}[sec5]{Lemma}
\begin{lem Gsyaei}
\label{lem Gsyaei}
Let $\l_1$, $\l_2 \in \bbZ$. 

(a) For all $r \in \bbN$, $g_1, \cdots, g_r \in Y_{1,2}$, we have
\[
 p_{\bf0}(g_{1} \cdots g_{r}v_{1,2}^{\bf0})=0
\q \rm{in $V_{1,2}(\l_1,\l_2)$}.
\]

(b) For all $r \in \bbN$, $i_1, \cdots, i_r \in \{1, 2\}$, we have
\[
 p_{\bf0}(f_{i_1,2} \cdots f_{i_r,2}v_{1,2}^{\bf0})=0
\q \rm{in $V_{1,2}(\l_1,\l_2)$}.
\]
\end{lem Gsyaei}

{\sl Proof.}
If we can prove (a), then we obtain (b) by 
(\ref{fact SHGf1}), (\ref{fact SHGf2}).
So we shall prove (a). 

Now we fix $r \in \bbN$, $g_1, \cdots, g_r \in Y_{1,2}$ 
and set $g:=g_1 \cdots g_r$. 
For $y \in Y_{1,2}$, we set 
\begin{eqnarray*} 
&&s(y):=\#\{1 \leq i \leq r \, | \, g_i=y\} \geq 0, 
\q m_g:=\sum_{i=1}^5 s(y_{i,2})\ep_{i,2}+s(y_{1,1})\ep_{1,1}, \\
&&W_g:=\bigoplus_{s=1}^r \bbC (g_sg_{s+1} \cdots g_rv_{1,2}^{\bf0}) 
\subset V_{1,2}(\l_1, \l_2).
\end{eqnarray*}
Then, $g v_{1,2}^{\bf0} \in \bbC v_{1,2}(m_g)$ 
by (\ref{def SHGy}), (\ref{fact SHGf1}), (\ref{fact SHGf2}). 
Since $\sum_{i=1}^5s(y_{i,2})+s(y_{1,1})=r>0$, 
there exists a $1 \leq i \leq 5$ such that $s(y_{i,2})>0$ 
or $s(y_{1,1})>0$. 

Case 1) $s(y_{1,1})>0 $: 
For $1 \leq r^{'} \leq r$, let $m^{(r^{'})} \in \bbZ_l^6$ such that 
$g_s g_{s+1} \cdots g_r \in \bbC v_{1,2}(m^{(r^{'})})$.
Let $1 \leq r_1 \leq r$ such that 
$g_{r_1}=y_{1,1}$ and $g_{r_1+1}, \cdots , g_{r} \neq y_{1,1}$. 
Then, by (\ref{def SHGy}), $m^{(r_1+1)}_{1,1}=0$. 
Hence, by the definition of $y_{1,1}$ in (\ref{def SHGy}), we get 
\begin{eqnarray*}
 g_{r_1}g_{r_1+1} \cdots g_r v_{1,2}^{\bf0} 
\in \bbC[-\l_1]v_{1,2}(m^{(r_1+1)}+\ep_{1,1}).
\end{eqnarray*}
Similarly, for $1 \leq r_2 < r_1$ such that 
$g_{r_2}=y_{1,1}$ and $g_{r_2+1}, \cdots , g_{r_1-1} \neq y_{1,1}$, 
we have 
\begin{eqnarray*}
 g_{r_2}g_{r_2+1} \cdots g_r v_{1,2}^{\bf0} 
\in \bbC[-\l_1+1][-\l_1]v_{1,2}(m^{(r_2+1)}+2\ep_{1,1}).
\end{eqnarray*}
By repeating this, we obtain 
\begin{eqnarray*} 
 g v_{1,2}^{\bf0} \in 
\bbC [-\l_1+s(y_{1,1})-1] \cdots [-\l_1+1][-\l_1]v_{1,2}(m_g).
\end{eqnarray*}  
Since $\l_1 \in \bbZ$ and $[l]=0$, 
if $s(y_{1,1}) \geq l$, then 
$[-\l_1+s(y_{1,1}-1] \cdots [-\l_1+1][-\l_1]=0$.
On the other hand, if $0<s(y_{1,1}) <l$, 
then $p_{\bf0}(v_{1,2}(m_g))=0$. 
Therefore, we obtain $p_{\bf0}(g v_{1,2}^{\bf0})=0$. 

Case 2) $s(y_{1,1})=0$ and $s(y_{5,2}) >0$: 
Since $s(y_{1,1})=0$, for all $1 \leq r^{'} \leq r$, $m^{(r^{'})}_{1,1}=0$.
Hence, we get
\begin{eqnarray*} 
y_{5,2}v_{1,2}(m_5, m_{1,1})
=[m_{5,2}-\l_2]_{\e_2}v_{1,2}(m_5+\ep_{5,2},m_{1,1}) 
\q \rm{in $W_g$}. 
\end{eqnarray*}
Thus, by the similar way to the proof of Case 1, 
we obtain $p_{\bf0}(g v_{1,2}^{\bf0})=0$. 

Case 3) There exists a $1 \leq i \leq 4$ such that 
$s(y_{1,1})=s(y_{5,2})= \cdots =s(y_{i+1,2})=0$ and $s(y_{i,2}) >0$: 
In this case , for all $1 \leq r^{'} \leq r$, 
$m_{1,1}^{(r^{'})}=m_{5,2}^{(r^{'})}= \cdots =m_{i+1,2}^{(r^{'})}=0$. 
Hence we have  
\begin{eqnarray*} 
y_{i,2}v_{1,2}(m_5, m_{1,1})
=[m_{i,2}-\l_{\tilde{i}}]_{\e_{\tilde{i}}}v_{1,2}(m_5+\ep_{i,2},m_{1,1}) 
\q \rm{in $W_g$}, 
\end{eqnarray*}
where $\tilde{i}:=1$ if $i=2,4$ and $\tilde{i}:=2$ if $i=1,3$. 
Therefore, by the similar way to the proof of Case 1, 
we obtain $p_{\bf0}(g v_{1,2}^{\bf0})=0$. 

By Case1--3, we obtain $p_{\bf0}(gv_{1,2}^{\bf0})=0$.
\qed 
\newtheorem{lemma Gf}[sec5]{Lemma}
\begin{lemma Gf}
\label{lemma Gf}
For all $\l_1$, $\l_2 \in \bbZ$, $\al \in \Delta_{+}$, we have 
$f_{\al,2}^lv_{1,2}^{\bf0}=0$ in $V_{k,n}(\l)$.
\end{lemma Gf} 

{\sl Proof.}
By Lemma \ref{lem Gsyaei}(b) and Proposition \ref{pro deg}, 
we obtain $p_{\bf0}(f_{\al,2}^lv_{1,2}^{\bf0})=0$. 
On the other hand, 
by Proposition \ref{pro CE}, \ref{sec5},
\[
 e_{i,2}f_{\al,2}^lv_{1,2}^{\bf0}=f_{\al,2}^le_{i,2}v_{1,2}^{\bf0}=0 
\q (i=1,2).
\]
Hence, by Proposition \ref{sec5}, we get
$f_{\al,2}^lv_{1,2}^{\bf0} \in \bbC v_{1,2}^{\bf0}$. 
Therefore, we obtain 
\begin{eqnarray*}
f_{\al,2}^lv_{1,2}^{\bf0}  
=p_{\bf0}(f_{\al,2}^lv_{1,2}^{\bf0})=0.
\end{eqnarray*}
\qed \\ 

Now, we construct nilpotent $\ue(G_2)$-modules (see \S 3).
For $\l_1$, $\l_2 \in \bbC$, 
let $L_{1,2}(\l_1, \l_2)$ be the $\ue(G_2)$-submodule of 
$V_{1,2}(\l_1, \l_2)$ generated by $v_{1,2}^{\bf0}$.
\newtheorem{thm nilG}[sec5]{Theorem}
\begin{thm nilG}
\label{thm nilG}
For any $\l_1$, $\l_2 \in \bbZ_l$, 
$L_{1,2}(\l_1, \l_2)$ is isomorphic to 
$L_{\e}^{\rm{nil}}(\l_1, \l_2)$ as $\ue(G_2)$-module. 
\end{thm nilG}

{\sl Proof.} 
By Proposition \ref{sec5}, 
$e_{1,2}v_{1,2}^{\bf0}=e_{2,2}v_{1,2}^{\bf0}=0$. 
Moreover, by (\ref{fact SHGt1}), (\ref{fact SHGt2}), 
\[
 t_{i,2}v_{1,2}^{\bf0}=\e_i^{\l_i}v_{1,2}^{\bf0} \q (i=1,2).
\]
So $L_{1,2}(\l_1, \l_2)$ is a finite dimensional 
highest weight $\ue(G_2)$-module with highest weight $(\l_1, \l_2)$. 
On the other hand, 
by Lemma \ref{lemma Gf}, 
$f_{\al,2}^lv_{1,2}^{\bf0}=0$ for all $\al \in \Delta_{+}$. 
Moreover, by Proposition \ref{pro deg}, \ref{sec5}, we have
$e_{\al,2}^lv_{1,2}^{\bf0}=0$ for all $\al \in \Delta_{+}$. 
Hence, by Proposition \ref{pro CE}, 
$e_{\al,2}^l=f_{\al,2}^l=0$ on $L_{1,2}(\l_1, \l_2)$ 
for all $\al \in \Delta_{+}$.  
Thus $L_{1,2}(\l_1, \l_2)$ is a nilpotent $\ue(G_2)$-module. 
Therefore, by Proposition \ref{sec5} and Proposition \ref{pro AN}, 
we obtain this theorem. 
\qed \\ 

If $\l_1=0$, 
then we can construct $L_{\e}^{\rm{nil}}(\l_1, \l_2)$ more easily. 
For $\l \in \bbC$, let $\rho^G(\l)$ as in (\ref{def rhoGll}). 
For $m \in \bbZ_{+}$, 
let $(\pi_m, \bbC)$ be the trivial representation of $U_{\e}(\g_m)$,
that is, 
\begin{eqnarray}
 \pi_{m}(e_{i,m})=\pi_{m}(f_{i,m})=0, \q \pi_{m}(t_{i,m})=1 
\q (1 \leq i \leq m), 
\label{def pi}
\end{eqnarray}
where $e_{i,0}:=f_{i,0}:=0$, $t_{i,0}:=1$, $\ue(\g_0):=\bbC$. 
For $\l_2 \in \bbC$, we define  
\begin{eqnarray*}
\phi_{2,2}:=\phi_{2,2}(\l_2):=\pi_1 \circ \rho^G(\nu_2^{\l_2})
:U_{\e}(G_2) \longrightarrow \rm{End}(V_5),
\end{eqnarray*}
where $\nu_{2}^{\l_2}:=3\l_2+9$.
We denote the $\ue(G_2)$-module 
associated with $(\phi_{2,2}(\l_2), V_5)$ by $V_{2,2}(\l_2)$. 
Let $L_{2,2}(\l_2)$ be the $\ue(G_2)$-submodule of $V_{2,2}(\l_2)$ 
generated by $v_{2,2}^{\bf0}:=v_5(0, \cdots, 0)$. 
Then, by the similar way to the proof of 
Proposition \ref{sec5}, Lemma \ref{lem Gsyaei}, \ref{lemma Gf}, 
and Theorem \ref{thm nilG}, 
we obtain the following proposition.
\newtheorem{pro nilG}[sec5]{Proposition}
\begin{pro nilG}
\label{pro nilG}
For any $\l_2 \in \bbZ_l$, 
$L_{2,2}(\l_2)$ is isomorphic to $L_l^{nil}(0, \l_2)$ 
as $\ue(G_2)$-module. 
In particular, 
for any $\l_2 \in \bbC$, 
we have $P(V_{2,2}(\l_2))=\bbC v_{2,2}^{\bf0}$. 
\end{pro nilG}
\section{Inductive construction of $\Lnil (\l)$ (Type B-case)}

\setcounter{equation}{0}
\renewcommand{\theequation}{\thesection.\arabic{equation}}

In this section, 
we construct all finite dimensional irreducible nilpotent 
$\ue(B_n)$-modules of type 1 inductively 
by using the Schnizer-homomorphisms of Theorem \ref{thm SH}(b). 

We set 
$a_n^{(0)}=(a_{i,n}^{(0)})_{i=1}^n$, 
$\ta_{n-1}^{(0)}=(\ta_{i,n-1}^{(0)})_{i=1}^{n-1}$, 
$b_n^{(0)}=(b_{i,n}^{(0)})_{i=1}^n$, 
$\tb_{n-1}^{(0)}=(\tb_{i,n-1}^{(0)})_{i=1}^{n-1} 
\in \bbC^{n}$ by 
\begin{eqnarray}
 a_{i,n}^{(0)}:=\ta_{i,n}^{(0)}:=1, 
\q  b_{i,n}^{(0)}:=n-i+1 \, (i \neq 1), 
\q b_{1,n}^{(0)}:=2n-1, 
 \q \tb_{i,n-1}^{(0)}:=i+n-2.
\label{sign ab}
\end{eqnarray} 
We fix $k \in I$. 
For $\l=(\l_{k}, \cdots ,\l_{n}) \in \bbC^{n-k+1}$, 
we define $\nu^{\l}=(\nu_{k}^{\l}, \cdots , \nu_{n}^{\l}) 
\in \bbC^{n-k+1}$ by 
\begin{eqnarray*}
 &&\nu_{i}^{\l}:=-2i+1-\sum_{j=k}^i\l_j \q (k \geq 2),  
\q \nu_{i}^{\l}
:=-2i+1-\frac{1}{2}\l_1-\sum_{j=2}^i\l_j \q (k=1),  
\end{eqnarray*}
where $k \leq i \leq n$. 
For $\l \in \bbC$, 
we define $\rho_{n}^{B}(\l)
:=\rho_n^B(a_n^{(0)}, \ta_{n-1}^{(0)},b_n^{(0)}, \tb_{n-1}^{(0)}, \l) 
: U_{\e}(B_n) \longrightarrow \rm{End}(V_n \otimes \tV_{n-1}) \otimes
U_{\e}(B_{n-1})$ 
(see Theorem \ref{thm SH}(b)),  
and let $(\pi_{k-1}, \bbC)$ be as in (\ref{def pi}).
We set  
\[
 V_{k,n}:=\bigotimes_{j=k}^n(V_j \otimes \tV_{j-1}). 
\] 
For $ \l=(\l_k, \cdots, \l_n) \in \bbC^{n-k+1}$, 
we define a $\ub$-representation 
$\phi_{k,n}:=\phi_{k,n}(\l):
\ub \longrightarrow  \textrm{End} (V_{k,n})$ by
\begin{equation}
\phi_{k,n}(\l):= 
\pi_{k-1} \circ  \rho_{k}^{B}(\nu_k^{\l})
 \circ \cdots  \circ  \rho_{n}^{B}(\nu_n^{\l}), 
\label{sig phikn}
\end{equation}
and denote the $\ue(B_n)$-module associated with 
$(\phi_{k,n}(\l), V_{k,n})$ by $V_{k,n}(\l)$. 

Let $m_n=(m_{1,n}, \cdots, m_{n,n}) \in \bbZ_l^n$, 
$\tm_{n-1}=(\tm_{1,n-1}, \cdots , \tm_{n-1,n-1}) \in \bbZ_l^{n-1}$, 
$w \in V_{k,n-1}$, 
$v=v_n(m_n) \otimes \tv_{n-1}(\tm_{n-1}) \otimes w 
\in V_{k,n}(\l)$. 
Then, by (\ref{def xz}), (\ref{def z}), (\ref{thm SHB1}), 
(\ref{sign ab}), for any $1<i<n$,
we have 
\begin{eqnarray}
e_{n,n}v
  &=&[-m_{n,n}](v_n(m_n-\ep_{n,n}) 
  \otimes \tv_{n-1}(\tm_{n-1}) \otimes w),
  \label{fact SHBen} \\ 
e_{i,n}v
  &=&[m_{i+1,n}-m_{i,n}]
  (v_n(m_n-\ep_{i,n}) \otimes \tv_{n-1}(\tm_{n-1}) \otimes w) \no \\
&&+[\tm_{i-1,n-1}-\tm_{i,n-1}]
  (v_n(m_n+\ep_{i+1,n}-\ep_{i,n}) \otimes 
  \tv_{n-1}(\tm_{n-1}-\tep_{i,n-1}) \otimes w)  \no\\
&&+v_n(m_n+\ep_{i+1,n}-\ep_{i,n}) \otimes 
  \tv_{n-1}(\tm_{n-1}-\tep_{i,n-1}+\tep_{i-1,n-1})\otimes (e_{i,n-1}w),
  \label{fact SHBei} \\
e_{1,n}v&=&[2m_{2,n}-m_{1,n}]_{\e_1}
  (v_n(m_n-\ep_{1,n}) \otimes \tv_{n-1}(\tm_{n-1}) \otimes w) \no \\
&&+[m_{1,n}-2\tm_{1,n-1}]_{\e_1}
  v_n(m_n+\ep_{2,n}-\ep_{1,n}) \otimes 
  \tv_{n-1}(\tm_{n-1}-\tep_{1,n-1}) \otimes w  \no \\
&&+v_n(m_n+\ep_{2,n})
 \otimes \tv_{n-1}(\tm_{n-1}-\tep_{1,n-1}) \otimes (e_{1,n-1}w). 
  \label{fact SHBe1} 
\end{eqnarray}

Let $m=(m_{i,j})_{1 \leq i \leq j, k \leq j \leq n } 
\in \bbZ_l^{N_n-N_{k-1}}$, 
$\tm=(\tm_{i,j-1})_{1 \leq i \leq j-1, k \leq j \leq n} 
\in \bbZ_l^{N_{n-1}-N_{k-2}}$, 
where $N_i:=\frac{1}{2}i(i+1)$ for $i \in \bbN$.
We set 
\begin{eqnarray}
&&v_{k,n}(m, \tm):=(\bigotimes_{j=k}^n v_j(m_{1,j}, \cdots m_{j,j}))
\otimes (\bigotimes_{j=k}^n \tv_{j-1}(\tm_{1,j-1}, 
\cdots \tm_{j-1,j-1})), \no \\
&&v_{k,n}^{\bf0}:=v_{k,n}(\bf0,\bf0).
\label{vkn} 
\end{eqnarray}
Let $P(V_{k,n}(\l))$ as in Definition \ref{def HWM} (i).
\newtheorem{pro pv}{Proposition}[section]
\begin{pro pv}
\label{pro pv}
For all $\l \in \bbC^{n-k+1}$, 
we obtain $P(V_{k,n}(\l))=\bbC v_{k,n}^{\bf0}$.
\end{pro pv}

{\sl Proof.} 
Since the actions of $e_{i,n}$ on $V_{k,n}(\l)$ do not depend on $\l$, 
we simply denote $V_{k,n}(\l)$ by $V_{k,n}$. 
By (\ref{fact SHBei}), (\ref{fact SHBe1}), obviously, 
$\bbC v_{k,n}^{\bf0} \subset P(V_{k,n})$. 
So we shall prove $P(V_{k,n}) \subset \bbC v_{k,n}^{\bf0}$
by induction on $n$. 

We assume $n=1$. Then we have $k=1$. 
Let $v=\sum_{m_{1,1} \in \bbZ_l}c(m_{1,1})v(m_{1,1}) \in V_{1,1}$ 
$(c(m_{1,1}) \in \bbC)$, 
and we assume $e_{1,1}v=0$. 
Then, by (\ref{fact SHBe1}), we get 
\begin{eqnarray*}
 0=e_{1,1}v 
=\sum_{m_{1,1} \in \bbZ_l}c(m_{1,1})[-m_{1,1}]_{\e_1}v_{1,1}(m_{1,1}-1). 
\end{eqnarray*}  
Hence, we obtain $c(m_{1,1})=0$ if $m_{1,1} \neq 0$. 
Therefore we have 
$v=c(0)v_{1,1}(0) \in \bbC v_{1,1}(0)=\bbC v_{1,1}^{\bf0}$. 

Now, we assume that $n >1$ 
and we obtain the case of $(n-1)$. 
Let  
\begin{eqnarray*}
 v=\sum_{m_n \in \bbZ_l^n, \tm_{n-1} \in \bbZ_l^{n-1}}
c(m_n,\tm_{n-1})(v_n(m_n) \otimes \tv_{n-1}(\tm_{n-1}) 
\otimes v_{m_n,\tm_{n-1}})\in V_{k,n},
\end{eqnarray*}
where 
$c(m_n,\tm_{n-1}) \in \bbC$, $v_{m_n,\tm_{n-1}} \in V_{k,n-1}$ 
($V_{n,n-1}:=\bbC v_{n,n-1}^{\bf0}$, $v_{n,n-1}^{\bf0}:=1$). 
We assume that $e_{i,n}v=0$ for all $1 \leq i \leq n$. 

First, we shall prove that 
$c(m_n,\tm_{n-1})=0$ if $m_n \neq \bf0$.
By (\ref{fact SHBen}), we get
\begin{eqnarray*} 
&&0=e_{n,n}v=
\sum_{m_n, \tm_{n-1}}c(m_n,\tm_{n-1})
[-m_{n,n}](v_n(m_n-\ep_{n,n}) \otimes \tv_{n-1}(\tm_{n-1}) \otimes 
v_{m_n,\tm_{n-1}}). 
\end{eqnarray*}
Hence, we obtain $c(m_n, \tm_{n-1})=0$ if $m_{n,n} \neq 0$. 
Now, we assume that there exists a $2 \leq i \leq n-1$ 
such that $c(m_n,\tm_{n-1})=0$ 
if $m_{i+1,n} \neq 0, \cdots m_{n-1,n} \neq 0$ or $m_{n,n} \neq 0$.
Then, by (\ref{fact SHBei}), we have
\begin{eqnarray*} 
&& 0=e_{i,n}v
=\sum_{m_n,\tm_{n-1}}c(m_n,\tm_{n-1})\{
[-m_{i,n}](v_n(m_n-\ep_{i,n}) \otimes \tv_{n-1}(\tm_{n-1}) \otimes
v_{m_n,\tm_{n-1}}) \\
&&\q +[\tm_{i-1,n-1}-\tm_{i,n-1}]
(v_n(m_n+\ep_{i+1,n}-\ep_{i,n}) \otimes 
\tv_{n-1}(\tm_{n-1}-\tep_{i,n-1}) \otimes v_{m_n,\tm_{n-1}}) \\
&&\q +v_n(m_n+\ep_{i+1,n}-\ep_{i,n}) \otimes 
\tv_{n-1}(\tm_{n-1}-\tep_{i,n-1}+\tep_{i-1,n-1}) \otimes 
(e_{i,n-1}v_{m_n,\tm_{n-1}}) \}.
\end{eqnarray*}
If $m_{i+1,n}=0$, then the $(i+1,n)$-component of 
$(m_n-\ep_{i,n})$ is $0$, 
and the one of $(m_n+\ep_{i+1,n}-\ep_{i,n})$ is $1$. 
Thus, by the linearly independence, 
$c(m_n,\tm_{n-1})=0$ if $m_{i,n} \neq 0$. 
Therefore we obtain 
$c(m_n,\tm_{n-1})=0$ 
if $m_{2,n} \neq 0, \cdots m_{n-1,n} \neq 0$ or $m_{n,n} \neq 0$ 
inductively.
Similarly, we have $c(m_n,\tm_{n-1})=0$ if $m_{1,n} \neq 0$ 
by using $e_{1,n}v=0$. 
Hence, we obtain $c(m_n,\tm_{n-1})=0$ if $m_n \neq \bf0$.
Therefore we get 
\begin{eqnarray*} 
 v=\sum_{ \tm_{n-1} \in \bbZ_l^{n-1}}
c(\bf0,\tm_{n-1})(v_n(\bf0) \otimes \tv_{n-1}(\tm_{n-1}) 
\otimes v_{\bf0,\tm_{n-1}}).
\end{eqnarray*}
Moreover, we have
\begin{eqnarray*} 
&&0=e_{i,n}v \\
&&=\sum_{ \tm_{n-1}}c(\bf0,\tm_{n-1})
\{ [\tm_{i-1,n-1}-\tm_{i,n-1}]
(v_n(\ep_{i+1,n}-\ep_{i,n}) \otimes 
\tv_{n-1}(\tm_{n-1}-\tep_{i,n-1}) \otimes v_{\bf0,\tm_{n-1}}) \\
&&\qq +v_n(\ep_{i+1,n}-\ep_{i,n}) \otimes 
\tv_{n-1}(\tm_{n-1}-\tep_{i,n-1}+\tep_{i-1,n-1}) \otimes 
(e_{i,n-1}v_{\bf0,\tm_{n-1}})\} \, (i \neq1),\\
&&0=e_{1,n}v \\
&&=\sum_{ \tm_{n-1}}c(\bf0,\tm_{n-1})
\{[-2\tm_{1,n-1}]_{\e_1}
(v_n(\ep_{2,n}-\ep_{1,n}) \otimes 
\tv_{n-1}(\tm_{n-1}-\tep_{1,n-1}) \otimes v_{\bf0,\tm_{n-1}}) \\
&&\qq +v_n(\ep_{2,n}) \otimes \tv_{n-1}(\tm_{n-1}-\tep_{1,n-1}) \otimes 
(e_{1,n-1}v_{\bf0,\tm_{n-1}})\}.
\end{eqnarray*}

Then, by the similar manner to the above proof, 
we obtain $c(\bf0,\tm_{n-1})=0$ if $\tm_{n-1} \neq \bf0$ 
by using $e_{i,n}v=0$ $(1 \leq i \leq n)$. 
Finally, we get
\begin{eqnarray*} 
 &&v=c(\bf0,\bf0)(v_n(\bf0) \otimes \tv_{n-1}(\bf0) \otimes v_{\bf0,\bf0}), \\
&&0=e_{i,n}v=c(\bf0,\bf0)v_n(\ep_{i+1,n}-\ep_{i,n}) \otimes 
\tv_{n-1}(-\tep_{i,n-1}+\tep_{i-1,n-1}) \otimes 
(e_{i,n-1}v_{\bf0,\bf0}) \, (i \neq 1,n), \\
&&0= e_{1,n}v
=c(\bf0,\bf0)v_n(\ep_{2,n}) \otimes \tv_{n-1}(-\tep_{1,n-1}) \otimes 
(e_{1,n-1}v_{\bf0,\bf0}).
\end{eqnarray*}
Hence $e_{i,n-1}v_{\bf0,\bf0}=0$ in $V_{k,n-1}$ 
for all $1 \leq i \leq n-1$ if $c(\bf0,\bf0) \neq 0$. 
So, by the assumption of the induction on $n$, 
we obtain $v_{\bf0, \bf0} \in \bbC v_{k,n-1}^{\bf0}$ 
if $c(\bf0,\bf0) \neq 0$. 
Therefore 
\begin{eqnarray*} 
 v\in \bbC (v_n(\bf0) \otimes \tv_{n-1}(\bf0) \otimes v_{k,n-1}^{\bf0})
=\bbC v_{k.n}^{\bf0}. 
\end{eqnarray*}
\qed \\ 

Let $m=(m_{i,j})_{1 \leq i \leq j, k \leq j \leq n } 
\in \bbZ_l^{N_n-N_{k-1}}$, 
$\tm=(\tm_{i,j-1})_{1 \leq i \leq j-1, k \leq j \leq n} 
\in \bbZ_l^{N_{n-1}-N_{k-2}}$. 
For $1 \leq i \leq  n$, $\rm{max}(k,i) \leq j \leq n$, 
we define
\begin{eqnarray}
&&\nu_{i,j}(m,\tm)
  :=m_{i+1,j}-2m_{i,j}+m_{i-1,j}
  +\tm_{i-2,j-1}-2\tm_{i-1,j-1}+\tm_{i,j-1} 
  \q (i \neq 1), \no \\
&&\nu_{1,j}(m,\tm)
  :=m_{2,j}-m_{1,j}+\tm_{1,j-1}, \no \\
&& \mu_{i,i-1}(m,\tm)
  :=\xi(i>k)(m_{i-1,i-1}+\tm_{i-2,i-2}), \no \\
&& \mu_{i,j}(m,\tm)
  :=\mu_{i,i-1}(m,\tm)+\sum_{r=\rm{max}(k,i)}^{j} \nu_{i,r}(m,\tm),
\label{fact SHBnu}
\end{eqnarray}
where
\begin{eqnarray}
 \xi(i>j):=
\begin{cases}
1 & (i>j) \\
0 & (i \leq j), 
\end{cases}
\q 
 \xi(i \geq j):=
\begin{cases}
1 & (i \geq j) \\
0 & (i < j). 
\end{cases}
\label{sig xi}
\end{eqnarray}
Let $v_{k,n}(m, \tm)$ be as in (\ref{vkn}). 
Then, by (\ref{def xz}), (\ref{thm SHB2}), 
(\ref{fact tB}), (\ref{sign ab}), 
we obtain 
\begin{eqnarray}
t_{i,n}v_{k,n}(m,\tm) 
&=&\e_i^{\mu_{i,n}(m, \tm)+\xi(i \geq k)\l_i}v_{k,n}(m,\tm).
\label{fact SHBt} 
\end{eqnarray}
In particular, we obtain the following lemma. 
\newtheorem{lem HW}[pro pv]{Lemma}
\begin{lem HW}
\label{lem HW}
For any $\l=(\l_k, \cdots , \l_n) \in \bbC^{n-k+1}$, 
$i \in I$, we obtain 
\begin{eqnarray*}
 t_{i,n}v_{k,n}^{\bf0}=\e_i^{\xi(i \geq k)\l_i}v_{k,n}^{\bf0} 
\q \rm{in $V_{k,n}(\l)$}.
\end{eqnarray*}
\end{lem HW}

Let $m=(m_{i,j})_{1 \leq i \leq j, k \leq j \leq n } 
\in \bbZ_l^{N_n-N_{k-1}}$, 
$\tm=(\tm_{i,j-1})_{1 \leq i \leq j-1, k \leq j \leq n} 
\in \bbZ_l^{N_{n-1}-N_{k-2}}$, 
$\l=(\l_k, \cdots, \l_n) \in \bbC^{n-k+1}$.  
For $i \in I$, $\rm{max}(k,i) \leq j \leq n$, 
we define  $y_{i,j}$, $\ty_{i-1,j-1} \in \rm{End}(V_{k,n}(\l))$ by: 
for $v=v_{k,n}(m, \tm)$,
\begin{eqnarray}
\ty_{i-1,j-1}v
  &:=&[\tm_{i-1,j-1}-\tm_{i,j-1}-\mu_{i,j-1}(m,\tm)-\xi(i \geq k)\l_i]
  v_{k,n}(m,\tm+\tep_{i-1,j-1}), \no \\
y_{i,j}v
  &:=&[m_{i+1,j}-m_{i,j}-\mu_{i,j}(m,\tm)-\xi(i \geq k)\l_i]
  v_{k,n}(m+\ep_{i,j},\tm)  \q (i \neq 1),\no  \\
y_{1,j}v
  &:=&[m_{1,j}-2\tm_{1,j-1}-\mu_{1,j}(m,\tm)-\xi(1 \geq k)\l_i]_{\e_1}
  v_{k,n}(m+\ep_{1,j},\tm), 
\label{sign SHBy}
\end{eqnarray}  
where $\ty_{0,j-1}:=0$. We set 
\begin{eqnarray*} 
Y_{k,n}:=\{y_{i,j}, \ty_{i-1,j-1} \, 
| \, i \in I, \rm{max}(k,i) \leq j \leq n \}.
\end{eqnarray*}
Then, by (\ref{def xz}), (\ref{def z}), (\ref{thm SHB3}), 
(\ref{sign ab}), we have
\begin{eqnarray} 
f_{i,n}v_{k,n}(m,\tm)
=\sum_{j=\rm{max}(k,i)}^n 
(y_{i,j}+\ty_{i-1,j-1})v_{k,n}(m,\tm) 
\q \rm{in $V_{k,n}(\l)$}. 
\label{sign Y}
\end{eqnarray}
Let $p_{\bf0}: V_{k,n}(\l)\longrightarrow \bbC v_{k,n}^{\bf0}$ 
be the projection.  
\newtheorem{lem syaei}[pro pv]{Lemma}
\begin{lem syaei}
\label{lem syaei}
Let $\l=(\l_k, \cdots, \l_n) \in \bbZ^{n-k+1}$. 

(a) For all $r \in \bbN$, $g_1, \cdots, g_r \in Y_{k,n}$, we have
\[
 p_{\bf0}(g_{1} \cdots g_{r}v_{k,n}^{\bf0})=0
\q in \, V_{k,n}(\l).
\]

(b) For all $r \in \bbN$, $i_1, \cdots, i_r \in I$, we have
\[
 p_{\bf0}(f_{i_1,n} \cdots f_{i_r,n}v_{k,n}^{\bf0})=0
\q in \, V_{k,n}(\l).
\]
\end{lem syaei}

{\sl Proof.}
If we can prove (a), then we obtain (b) by (\ref{sign Y}). 
So we shall prove (a). 

Let $r \in \bbN$, $g_1, \cdots, g_r \in (Y_{k,n}-\{0\})$ 
and set $g:=g_1 \cdots g_r$. 
For $y \in Y_{k,n}$, we set 
\begin{eqnarray*} 
&&s(y):=\#\{1 \leq r^{'} \leq r \, | \, g_{r^{'}}=y\} \geq 0, 
\q s_j:=\sum_{i=1}^j (s(y_{i,j})+s(\ty_{i-1,j-1})) 
\q (k \leq j \leq n), \\
&&m_g:=\sum_{j=k}^n \sum_{i=1}^j 
(s(y_{i,j})\ep_{i,j}+ s(\ty_{i-1,j-1})\tep_{i-1,j-1}), 
\q W_g:=\bigoplus_{r^{'}=1}^r \bbC 
(g_{r^{'}}g_{r^{'}+1} \cdots g_rv_{k,n}^{\bf0}).
\end{eqnarray*}
Then, $g v_{k,n}^{\bf0} \in \bbC v_{k,n}(m_g)$ 
by (\ref{sign SHBy}). 
Since $\sum_{j=k}^n s_j=r >0$, 
there exists a $k \leq j \leq n$ such that 
$s_{k}= \cdots =s_{j-1}=0$ and $s_j>0$. 
Then, $s(y_{p,q})=s(\ty_{p-1,q-1})=0$ 
for all $k \leq q <j$, $1 \leq p \leq q$. 
Thus, for any $1 \leq  r^{'} \leq r$, 
there exist $m^{(r^{'})} \in \bbZ_l^{N_n-N_{k-1}}$, 
$\tm^{(r^{'})} \in \bbZ_l^{N_{n-1}-N_{k-2}}$ such that 
$g_{r^{'}} g_{r^{'}+1} \cdots g_rv_{k,n}^{\bf0} 
\in \bbC v_{k,n}(m^{(r^{'})}, \tm^{(r^{'})})$ 
and $m_{p,q}^{(r^{'})}=\tm_{p-1,q-1}^{(r^{'})}=0$ 
for all $k \leq q <j$, $1 \leq p \leq q$. 
Hence, by (\ref{fact SHBnu}), (\ref{sign SHBy}), in $W_g$,
\begin{eqnarray} 
y_{i,j}v_{k,n}(m,\tm)
&=&[m_{i+1,j}-m_{i,j}-\nu_{i,j}(m, \tm)-\xi(i \geq k)\l_i] 
v_{k,n}(m+\ep_{i,j},\tm), \no \\
y_{1,j}v_{k,n}(m,\tm)
&=&[m_{1,j}-2\tm_{1,j-1}-\xi(1 \geq k)\l_1]_{\e_1}
v_{k,n}(m+\ep_{1,j},\tm), \no \\
\ty_{i-1,j-1}v_{k,n}(m,\tm)
&=&[\tm_{i-1,j-1}-\tm_{i,j-1}-\xi(i \geq k)\l_i]
v_{k,n}(m,\tm+\tep_{i-1,j-1}),
\label{sig yonW} 
\end{eqnarray}
for $2 \leq i \leq j$, 
where $\xi(i \geq j)$ as in (\ref{sig xi}). 

On the other hand, since $s_j>0$, 
there exist $i$ $(1 \leq i \leq j)$ such that 
$s(y_{i,j})>0$ or $s(\ty_{i-1,j-1})>0$. 
Now, we assume $s(\ty_{j-1,j-1})>0$. 
Let $r_1$ $(1 \leq r_1 \leq r)$ such that $g_{r_1}=\ty_{j-1,j-1}$ 
and $g_{r_1+1}, \cdots , g_r \neq \ty_{j-1,j-1}$. 
Then, by (\ref{sign SHBy}),  $\tm_{j-1,j-1}^{(r_1+1)}=0$.
Hence, by (\ref{sig yonW}), 
\begin{eqnarray*} 
g_{r_1}g_{r_1+1} \cdots g_{r}v_{k,n}^{\bf0}
\in \bbC[-\l_j]v_{k,n}(m^{(r_1+1)}, \tm^{(r_1+1)}+\tep_{j-1,j-1}).
\end{eqnarray*}
Thus, by the similar way to the proof of Case 1 in Lemma \ref{lem Gsyaei}, 
we have
\begin{eqnarray*} 
 g v_{k,n}^{\bf0} \in 
\bbC [-\l_j+s(\ty_{j-1,j-1})-1] \cdots [-\l_j+1][-\l_j]v_{k,n}(m_g), 
\end{eqnarray*}  
and $p_{\bf0}(g v_{k,n}^{\bf0})=0$. Similarly, 
if there exists $i$ $(2 \leq i \leq j-1)$ such that 
$s(\ty_{j-1,j-1})= \cdots =s(\ty_{i,j-1})=0, s(\ty_{i-1,j-1})>0$, 
then we have 
\[
 \ty_{i-1,j-1}v_{k,n}(m,\tm)
=[\tm_{i-1,j-1}-\delta_{i \geq k}\l_i]v_{k,n}(m,\tm+\tep_{i-1,j-1}) 
\q \rm{in $W_g$}, 
\] 
and $p_{\bf0}(g v_{k,n}^{\bf0})=0$.
If $s(\ty_{1,j-1})= \cdots =s(\ty_{j-1,j-1})=0$, 
and there exists $i$ $(1 \leq i \leq  j)$ such that 
$s(y_{1,j})= \cdots =s(y_{i-1,j})=0, s(y_{i,j})>0$, 
then we obtain
\begin{eqnarray} 
&&y_{i,j}v_{k,n}(m,\tm)
=[m_{i,j}-\delta_{i \geq k}\l_i]_{\e_i} 
v_{k,n}(m+\ep_{i,j},\tm) \q \rm{in $W_g$},\no
\end{eqnarray}
and $p_{\bf0}(g v_{k,n}^{\bf0})=0$. 

Consequently, we obtain that
 $p_{\bf0}(gv_{k,n}^{\bf0})=0$ if $s_j>0$. 
It amount to $p_{\bf0}(gv_{k,n}^{\bf0})=0$.
\qed 
\newtheorem{lem f}[pro pv]{Lemma}
\begin{lem f}
\label{lem f}
For all $\l=(\l_k, \cdots, \l_n) \in \bbZ^{n-k+1}$, 
$\al \in \Delta_{+}$, we have 
\begin{eqnarray}
 f_{\al,n}^lv_{k,n}^{\bf0}=0 \q \rm{in $V_{k,n}(\l)$}. \no
\end{eqnarray}
\end{lem f} 

{\sl Proof.}
By Lemma \ref{lem syaei}(b) and Proposition \ref{pro deg}, 
we obtain $p_{\bf0}(f_{\al,n}^lv_{k,n}^{\bf0})=0$. 
On the other hand, 
by Proposition \ref{pro CE}, \ref{pro pv},
\[
 e_{i,n}f_{\al,n}^lv_{k,n}^{\bf0}
=f_{\al,n}^le_{i,n}v_{k,n}^{\bf0}=0, 
\]
for all $i \in I$. 
Hence, by Proposition \ref{pro pv}, we get
$f_{\al,n}^lv_{k,n}^{\bf0} \in \bbC v_{k,n}^{\bf0}$. 
Therefore 
\begin{eqnarray*}
f_{\al,n}^lv_{k,n}^{\bf0}  
=p_{\bf0}(f_{\al,n}^lv_{k,n}^{\bf0})=0.
\end{eqnarray*}
\qed \\ 

Now, we construct nilpotent $\ue(B_n)$-modules (see \S 3).
For $\l \in \bbC^{n-k+1}$, 
let $L_{k,n}(\l)$ be the $\ue(B_n)$-submodule of 
$V_{k,n}(\l)$ generated by $v_{k,n}^{\bf0}$.
\newtheorem{thm nilB}[pro pv]{Theorem}
\begin{thm nilB}
\label{thm nilB}
For any $k \in I$, $\l=(\l_k, \cdots, \l_n) \in \bbZ_l^{n-k+1}$,
$L_{k,n}(\l)$ is isomorphic to 
$L_{\e}^{\rm{nil}}$$(0, \cdots, 0, \l_k, \cdots, \l_n)$ 
as $\ue(B_n)$-module. 
\end{thm nilB}

{\sl Proof.} 
By Proposition \ref{pro pv} and Lemma \ref{lem HW}, 
$L_{k,n}(\l)$ is a highest weight $\ue(B_n)$-module 
with highest weight $(0, \cdots, 0, \l_k, \cdots \l_n)$. 
On the other hand, 
by Lemma \ref{lem f}, $f_{\al,n}^lv_{k,n}^{\bf0}=0$ 
for all $\al \in \Delta_{+}$. 
Moreover, by Proposition \ref{pro deg}, \ref{pro pv},
$e_{\al,n}^lv_{k,n}^{\bf0}=0$ for all $\al \in \Delta_{+}$. 
Hence, by Proposition \ref{pro CE}, 
$e_{\al,n}^l=f_{\al,n}^l=0$ on $L_{k,n}(\l)$ for all $\al \in \Delta_{+}$.  
Thus $L_{k,n}(\l)$ is a nilpotent $\ue(B_n)$-module. 
Therefore, by Proposition \ref{pro pv} and Proposition \ref{pro AN}, 
we obtain this theorem. 
\qed \\ 

In particular, if $k=1$ 
then we obtain all finite dimensional irreducible nilpotent 
$\ue(B_n)$-modules of type 1. 
\section{Other cases}

\setcounter{equation}{0}
\renewcommand{\theequation}{\thesection.\arabic{equation}}

In this section, we construct 
$\ue(\g_n)$-modules $\Lnil (\l)$ in the case of 
$\g_n=A_n$, $C_n$ or $D_n$ 
inductively by using the Schnizer-homomorphisms 
in Theorem \ref{thm SH}(a), (c), (d). 

We set  
$a_n^{(0)}=(a_{i,n}^{(0)})_{i=1}^n$, 
$\ta_n^{(0)}=(\ta_{i,n}^{(0)})_{i=1}^n$, 
$b_n^{(0)}=(b_{i,n}^{(0)})_{i=1}^n$,  
$\tb_n^{(0)}=(\tb_{i,n}^{(0)})_{i=1}^n 
\in \bbC^{n}$ by, 
\begin{eqnarray*}
 && a_{i,n}^{(0)}:=\ta_{i,n}^{(0)}:=1 \q 
(\rm{$\g_n=A_n$, $C_n$, $D_n$}), \\
&& b_{i,n}^{(0)}:=i  \q (\g_n=A_n),  
\q b_{i,n}^{(0)}:=n-i+1 \q (\g_n=C_n),   \\
&& b_{i,n}^{(0)}:=n-i+1 \q (i \neq 1),
\q b_{1,n}^{(0)}:=n-1 \q (n \neq 1), 
\q b_{1,1}^{(0)}:=1 \q (\g_n=D_n),   \\
&& \tb_{i,n}^{(0)}:=n+i-1 \q (\rm{$\g_n=C_n$, $D_n$}).
\end{eqnarray*} 
We fix $k \in I$. 
For $\l=(\l_{k}, \cdots ,\l_{n}) \in \bbC^{n-k+1}$, 
we define $\nu^{\l}=(\nu_{k}^{\l}, \cdots , \nu_{n}^{\l}) 
\in \bbC^{n-k+1}$ by 
\begin{eqnarray*}
 &&\nu_{i}^{\l}:=-i-1-\frac{1}{i}\sum_{j=k}^i(j \l_j)\q (\g_n=A_n),  
\q\nu_{i}^{\l}:=-2i-\sum_{j=k}^i\l_j   \q (\g_n=C_n),  \\
 &&\nu_{i}^{\l}:=-2i+2-\sum_{j=k}^i\l_j \q (k \geq 3), 
\q \nu_{i}^{\l}:=-2i+3-\frac{1}{2}\l_2-\sum_{j=3}^i\l_j \q (k=2), \\
&& \nu_{i}^{\l}:=-2i+1-\frac{1}{2}\l_1-\frac{1}{2}\l_2-\sum_{j=3}^i\l_j 
\q (k=1), \q (\g_n=D_n), 
\end{eqnarray*}
where $k \leq i \leq n$.
For $\l \in \bbC$, 
we set $\rho_{n}^{A}(\l):=\rho_n^A(a_n^{(0)}, b_{n}^{(0)} \l)$, 
$\rho_{n}^{C}(\l)
:=\rho_n^C(a_n^{(0)}, \ta_{n-1}^{(0)},b_n^{(0)}, \tb_{n-1}^{(0)}, \l)$, 
and $\rho_{n}^{D}(\l)
:=\rho_n^D(a_n^{(0)}, \ta_{n-2}^{(0)},b_n^{(0)}, \tb_{n-2}^{(0)}, \l)$ 
(see Theorem \ref{thm SH}(a), (c), (d)).  
We define  
\begin{eqnarray*}
 &&V_{k,n}:=\bigotimes_{j=k}^nV_j\q (\g_n=A_n),    
\q V_{k,n}:=\bigotimes_{j=k}^n(V_j \otimes \tV_{j-1})\q (\g_n=C_n), \\
&& V_{k,n}:=\bigotimes_{j=k}^n(V_j \otimes \tV_{j-2})\q (\g_n=D_n). 
\end{eqnarray*}
For $ \l=(\l_k, \cdots, \l_n) \in \bbC^{n-k+1}$, 
we define $\ue (\g_n)$-representations 
$\phi_{k,n}:=\phi_{k,n}(\l):
 \ue(\g_n) \longrightarrow  \textrm{End} (V_{k,n})$ by
\begin{eqnarray*}
&&\phi_{k,n}(\l):= 
\pi_{k-1} \circ  \rho_{k}^{\g}(\nu_k^{\l})
 \circ \cdots  \circ  \rho_{n}^{\g}(\nu_n^{\l}). 
\end{eqnarray*}
We denote the $\ue(\g_n)$-module associated with 
$(\phi_{k,n}(\l), V_{k,n})$ by $V_{k,n}(\l)$.
We set 
\begin{eqnarray*}
&& v_{k,n}^{\bf0}:=\bigotimes_{j=k}^n v_j(0, \cdots, 0) \q (\g_n=A_n),  \\  
&& v_{k,n}^{\bf0}:=\bigotimes_{j=k}^n 
(v_j(0, \cdots, 0) \otimes \tv_{j-1}(0, \cdots, 0)) \q (\g_n=C_n),  \\  
&& v_{k,n}^{\bf0}:=\bigotimes_{j=k}^n 
(v_j(0, \cdots, 0) \otimes \tv_{j-2}(0, \cdots, 0))\q (\g_n=D_n).
\end{eqnarray*}

For $\l \in \bbC^{n-k+1}$, 
let $L_{k,n}(\l)$ be the $\ue(\g_n)$-submodule of 
$V_{k,n}(\l)$ generated by $v_{k,n}^{\bf0}$. 
Then we have the following theorem 
by the similar way to \S6.
\newtheorem{thm nil}{Theorem}[section]
\begin{thm nil}
\label{thm nil}
Let $\g_n=A_n$, $C_n$ or $D_n$.
Then, for any $k \in I$, $\l=(\l_k, \cdots, \l_n) \in \bbZ_l^{n-k+1}$, 
$L_{k,n}(\l)$ is isomorphic to $L_l^{nil}(0, \cdots, 0, \l_k, \cdots \l_n)$ 
as $\ue(\g_n)$-module. 
In particular, 
for any $\l \in \bbC^{n-k+1}$, 
we have $P(V_{k,n}(\l))=\bbC v_{k,n}^{\bf0}$. 
\end{thm nil}

Consequently,   
we obtain all finite dimensional irreducible nilpotent modules of type 1 
inductively by using the Schnizer homomorphisms 
for quantum algebras at roots of unity of 
type $A_n$, $B_n$, $C_n$, $D_n$ or $G_2$. \\

\textbf{Acknowledgement:}
I would like to thank T. Nakashima for helpful discussions.



\begin{thebibliography}{99}   
\bibitem{AN}
Y. Abe and T. Nakashima:
Nilpotent representations of classical quantum groups at roots of unity. 
J. Math. Phys. 46 (2005), No. 12, 113505 1-19.

\bibitem{CP} 
V. Chari and A. Pressley:
A Guide to Quantum Groups. 
Cambridge University Press, Cambridge. (1994).

\bibitem{DJMM} 
E. Date, M. Jimbo, K. Miki and T. Miwa: 
Cyclic Representations of $U_q(sl(n+1,\bbC))$ at $q^N=1$.
Publ. RIMS, Kyoto Univ. 27 (1991) 366-437.

\bibitem{DK} 
C. De Concini and V. G. Kac: 
Actes du Colloque en I'honneur de Jacques Diximier, 
edited by A.Connes, M Duflo, A. Joseph and R.Rentschler 
(Prog. Math. Birkhauser.). 
92 (1990) 471-506. 

\bibitem{J} J. C. Jantzen: 
Lectures on Quantum Groups. 
GSM. vol.6 (1996). 

\bibitem{L1}
G. Lusztig:  
Modular representations and quantum groups.
Contemp.Math 82 (1989) 59-77. 

\bibitem{L2}
G. Lusztig:
 Quantum groups at root of 1. 
Geom.Dedicata 35 (1990) 89-113. 

\bibitem{N}
T. Nakashima: 
Irreducible modules of finite dimensional quantum algebras 
of type A at roots of unity.
J.Math.Phys. vol.43 No.4 (2002) 2000-2014.

\bibitem{S1}
W. A. Schnizer: 
Roots of unity: Representations for symplectic and 
orthogonal quantum groups. 
 J.Math.Phys. 34 (1993) 4340-4363. 

\bibitem{S2}
W. A. Schnizer: 
Roots of unity: Representations of Quantum Group.
Commun. Math. Phys. 163 (1994) 293-306. 

\end{thebibliography}
\end{document}